
\input amstex
\documentstyle{amsppt}
\loadbold
\loadeufm
\loadeufb

\magnification=\magstep1
\hsize=6.5truein
\vsize=9truein

\document

\baselineskip=14pt

\font \smallrm=cmr10 at 10truept
\font \smallbf=cmbx10 at 10truept 
 at 10truept
\font \smallsl=cmsl10 at 10truept
\font\sc=cmcsc10

\def \loongrightarrow {\relbar\joinrel\relbar\joinrel\rightarrow}

\def \longtwoheadrightarrow
{\relbar\joinrel\relbar\joinrel\twoheadrightarrow}

\def \gerg {\frak g}
\def \gerb {\frak b}
\def \gert {\frak t}

\def \gergl {{\frak{gl}}}
\def \gersl {{\frak{sl}}}
\def \gergln {\gergl_{\hskip0,5pt n}}
\def \gersln {\gersl_{\hskip0,5pt n}}
\def \und1 {\underline{1}}

\def \Z {\Bbb{Z}}
\def \Q {\Bbb{Q}}
\def \keps {\Bbbk_\varepsilon}
\def \kq {\Bbbk(q)}
\def \kqqm {\Bbbk\big[q,q^{-1}\big]}
\def \qm {q^{-1}}
\def \uhat {\widehat{U}}
\def \utilde {\widetilde{U}}
\def \fhat {\widehat{F}}
\def \ftilde {\widetilde{F}}
\def \edot {\dot{E}}
\def \fdot {\dot{F}}
\def \uqg {U_q(\gerg)}
\def \uqmg {U_q^{\scriptscriptstyle M}\!(\gerg)}

\def \uqpgl {U_q^{\scriptscriptstyle P}\!(\gergln)}
\def \uqqgl {U_q^{\scriptscriptstyle Q}\!(\gergln)}

\def \uqpsl {U_q^{\scriptscriptstyle P}\!(\gersln)}
\def \uqqsl {U_q^{\scriptscriptstyle Q}\!(\gersln)}
\def \uqmt {U_q^{\scriptscriptstyle M}\!(\gert)}

\def \uqmbp {U_q^{\scriptscriptstyle M}\!(\gerb_+)}
\def \uqmpbp {U_q^{\scriptscriptstyle M'}\!(\gerb_+)}
\def \uqpbp {U_q^{\scriptscriptstyle P}\!(\gerb_+)}
\def \uqqbp {U_q^{\scriptscriptstyle Q}\!(\gerb_+)}
\def \uqmbm {U_q^{\scriptscriptstyle M}\!(\gerb_-)}
\def \uqmpbm {U_q^{\scriptscriptstyle M'}\!(\gerb_-)}
\def \uqpbm {U_q^{\scriptscriptstyle P}\!(\gerb_-)}
\def \uqqbm {U_q^{\scriptscriptstyle Q}\!(\gerb_-)}
\def \fqmg {F_q^{\scriptscriptstyle M}\!\big[G\big]}
\def \fqmpg {F_q^{\scriptscriptstyle M'}\!\big[G\big]}
\def \fqmbp {F_q^{\scriptscriptstyle M}\!\big[B_+\big]}
\def \fqmpbp {F_q^{\scriptscriptstyle M'}\!\big[B_+\big]}
\def \fqmbm {F_q^{\scriptscriptstyle M}\!\big[B_-\big]}
\def \fqmpbm {F_q^{\scriptscriptstyle M'}\!\big[B_-\big]}

\def \uqpbp {U_q^{\scriptscriptstyle P}\!(\gerb_+)}
\def \uqpbm {U_q^{\scriptscriptstyle P}\!(\gerb_-)}
\def \fqpg {F_q^{\scriptscriptstyle P}\!\big[G\big]}
\def \fqpbp {F_q^{\scriptscriptstyle P}\!\big[B_+\big]}
\def \fqpbm {F_q^{\scriptscriptstyle P}\!\big[B_-\big]}
\def \fqpgl {F_q^{\scriptscriptstyle P}\!\big[GL_n\big]}


\topmatter

{\ } 

\vskip-57pt  

\hfill   {\smallrm {\smallsl Proceedings of the Edinburgh
Mathematical Society\/}  {\smallbf 49}  (2006), 291--308}  
\hskip19pt   {\ }  

\vskip75pt  

\title
  PRESENTATION BY BOREL SUBALGEBRAS AND CHEVALLEY  \\   
  GENERATORS FOR QUANTUM ENVELOPING ALGEBRAS  
\endtitle

\author
       Fabio Gavarini
\endauthor

\rightheadtext{ CHEVALLEY GENERATORS FOR QUANTUM ENVELOPING ALGEBRAS }

\affil
  Universit\`a degli Studi di Roma ``Tor Vergata'' ---
Dipartimento di Matematica  \\
  Via della Ricerca Scientifica 1, I-00133 Roma --- ITALY  \\
\endaffil

\address\hskip-\parindent
  Fabio Gavarini  \newline
     \indent   Universit\`a degli Studi di Roma ``Tor Vergata''
---   Dipartimento di Matematica   \newline
     \indent   Via della Ricerca Scientifica 1,
I-00133 Roma, ITALY  \newline
     \indent   e-mail: \, gavarini\@{}mat.uniroma2.it
\endaddress

\abstract
   We provide an alternative approach to the
Faddeev-Reshetikhin-Takhtajan presentation of the
quantum group  $ \uqg $,  with  $ L $--operators
as generators and relations ruled by an  $ R $--matrix.
We look at  $ \uqg $  as being generated by the quantum Borel
subalgebras  $ U_q(\gerb_+) $  and  $ U_q(\gerb_-) $,  and use the
standard presentation of the latters as quantum function algebras.
When  $ \, \gerg = \gergl_n \, $  these Borel quantum function
algebras are generated by the entries of a triangular  $ q $--matrix,
thus eventually  $ U_q(\gergl_n) \, $  is generated by the entries
of an upper triangular and a lower triangular  $ q $--matrix, 
which share the same diagonal.  The same elements generate
over  $ \Bbbk\big[q,q^{-1}\big] $  the unrestricted  $ \Bbbk
\big[q,q^{-1}\big] $--integer  form of  $ U_q(\gergl_n) $  of
De Concini and Procesi, which we present explicitly, together
with a neat description of the associated quantum Frobenius
morphisms at roots of 1.  All this holds, {\sl mutatis mutandis}, 
for  $ \, \gerg = \gersl_n \, $  too.  
\endabstract

\endtopmatter

\footnote""{Keywords: \ {\sl Quantum Groups,  $ L $--operators,
Quantum Root Vectors}.}
\footnote""{ 2000 {\it Mathematics Subject Classification:} \
Primary 17B37, 20G42; Secondary 81R50. }

\vskip0,67truecm

\centerline {\bf Introduction }

\vskip10pt

   Let  $ \gerg $  be a semisimple Lie algebra over a field  $ \Bbbk
\, $.  Classically, it has two standard presentations: Serre's one,
which uses a minimal set of generators, and Chevalley's one, using a
linear basis as generating set.  If  $ \gerg $  instead is reductive
a presentation is obtained by that of its semisimple quotient by
adding the center.  When  $ \, \gerg = \gergln \, $,  Chevalley's
generators are the elementary matrices, and Serre's ones form a
distinguished subset of them; the general case of any classical
matrix Lie algebra  $ \gerg $  is a slight variation on this theme.
Finally, both presentations yield also presentations of  $ U(\gerg) $,
the universal enveloping algebra of  $ \gerg \, $.
                                               \par
   At the quantum level, one has correspondingly a Serre-like and a
Chevalley-like presentation of  $ \uqg $,  the quantized universal
enveloping algebra associated to  $ \gerg $  after Jimbo and Lusztig
(i.e.~defined over the field  $ \kq $,  where  $ q $  is an
indeterminate).  The first presentation is used by Jimbo (cf.~[Ji1])
and Lusztig (see [Lu2]) and,  {\sl mutatis mutandis},  by Drinfeld too;
in this case the generators are  $ q $--analogues  of the Serre's
generators, and starting from them one builds quantum root vectors
via two different methods: iterated quantum brackets, as in [Ji2]   ---
and maybe others ---   or braid group action, like in [Lu2]; see [Ga2]
for a comparison between these two methods.  The second presentation
was introduced by Faddeev, Reshetikhin and Takhtajan (in [FRT]): the
generators in this case, called  $ L $--operators,  are  $ q $--analogues
of the classical Chevalley generators; in particular, they are quantum
root vectors themselves.  An explicit comparison between quantum
Serre-like generators and  $ L $--operators  appears in [FRT], \S 2,
for the cases of  {\sl classical\/}  $ \gerg \, $;  \, on the other
hand, in [No], \S 1.2, a similar comparison is made for  $ \, \gerg
= \gergln \, $  between  $ L $--operators  and quantum root vectors
(for  {\sl any\/}  root) built out of Serre's generators.
%
%
                                               \par
   The first purpose of this note is to provide an alternative approach
to the FRT presentation of  $ \uqg \, $:  it amounts to a series of
elementary steps, yet the final outcome seems noteworthy.  As a
second, deeper result, we give an explicit presentation of the 
$ \kqqm $--subalgebra  of  $ \uqg \, $  generated by  $ L $--operators, 
call it  $ \utilde_q(\gerg) \, $.  By construction, this is nothing but
the  {\sl unrestricted\/}  $ \kqqm $--integer  form of  $ \uqg \, $, 
defined by De Concini and Procesi (see [DP]), whose semiclassical limit
is  $ \; \utilde_q(\gerg) \Big/ (q\!-\!1) \, \utilde_q(\gerg) \, \cong
\, F\big[G^*\big] \, $,  \; where  $ \, G^* \, $  is a connected Poisson
algebraic group dual to  $ \gerg \, $  (cf.~[DP], [Ga1] and [Ga3], \S
7.3 and \S 7.9): our explicit presentation of  $ \utilde_q(\gerg) $ 
yields another, independent (and much easier) proof of this fact. 
Third, by [DP] we know that  {\sl quantum Frobenius morphisms\/} 
exist, which embed  $ F\big[G^*\big] $  into the specializations
of  $ \utilde_q(\gerg) $  at roots of 1: then our presentation of 
$ \utilde_q(\gerg) $  provides an explicit description of them.   
                                               \par
   This analysis shows that the two presentations of  $ \uqg $  correspond
to different behaviors w.r.t.~to specializations.  Indeed, let  $ \uhat_q
(\gerg) $  be the  $ \kqqm $--algebra  given by Jimbo-Lusztig presentation 
{\sl over\/} $ \kqqm \, $.  Its specialization at  $ \, q = 1 \, $  is 
$ \; \uhat_q(\gerg) \Big/ (q\!-\!1) \, \uhat_q(\gerg) \, \cong \, U(\gerg)
\, $  (up to technicalities), with  $ \gerg $  inheriting a Lie bialgebra
structure (see [Ji1], [Lu2], [DL]).  On the other hand, the integer form 
$ \utilde_q(\gerg) $  mentioned above specializes to  $ F\big[G^*\big]
\, $,  the Poisson structure on  $ G^* $  being exactly the one dual to
the Lie bialgebra structure on  $ \gerg \, $.  So the existence of two
different presentations of  $ \uqg $  reflects the deep fact that  $ \uqg $ 
provides, taking suitable integer forms, quantizations of two different
semiclassical objects (this is a general fact, see [Ga3--4]).  To
the author's knowledge, this was not known so far, as the FRT
presentation of  $ \uqg $  was never used to study the integer
form  $ \utilde_q(\gerg) \, $.  
                                               \par
   Let's sketch in short the path we follow.  First, we note that 
$ \uqg $  is generated by the quantum Borel subgroups  $ U_q(\gerb_-) $ 
and  $ U_q(\gerb_+) $  (where  $ \gerb_- $  and  $ \gerb_+ $  are opposite
Borel subalgebras of  $ \gerg $),  which share a common copy of the
quantum Cartan subgroup  $ U_q(\frak{t}) \, $.  Second, there exist
Hopf algebra isomorphisms  $ \; U_q(\gerb_-) \cong F_q\big[B_-\big] \; $ 
and  $ \; U_q(\gerb_+) \cong F_q\big[B_+\big] \, $,  \; where  $ F_q
\big[B_-\big] $  and  $ F_q\big[B_+\big] $  are the quantum function
algebras associated to  $ \gerb_- $  and  $ \gerb_+ $  respectively. 
Third, when  $ \gerg $  is classical we resume the explicit presentation
by generators and relations of  $ F_q\big[B_-\big] $  and  $ F_q \big[
B_+ \big] $,  as given in [FRT], \S 1.  Fourth, from the above we argue
a presentation of  $ \uqg $  where the generators are those of  $ F_q
\big[B_-\big] $  and  $ F_q\big[B_+\big] $,  the toral ones being taken
only once, and relations are those of these quantum function algebras
plus some additional relations between generators of opposite quantum
Borel subgroups.  We perform this last step in full detail for  $ \,
\gerg = \gergln \, $  and, with slight changes, for  $ \, \gerg =
\gersln \, $  as well.  Fifth, we refine the last step to provide
a presentation of  $ \, \utilde_q(\gerg) \, $.
                                               \par
   As an application, our results apply also (with few changes)
to the Drinfeld-like quantum groups  $ U_\hbar(\gerg) \, $:  in
particular we get a presentation of an  $ \hbar $--deformation 
of  $ F[G^*] \, $,  say  $ \, \utilde_\hbar(\gerg) =: F_\hbar[G^*] \, $. 
An explicit gauge equivalence between this  $ F_\hbar[G^*] $  and the 
$ \hbar $--deformation  provided by Kontsevitch' recipe is given in [FG].  

\vskip15pt

\centerline{ \sc acknowledgements }

\centerline{ The author thanks P.~M\"{o}seneder Frajria and
D.~Parashar for helpful conversations. }

\vskip1,3truecm

\centerline {\bf \S \; 1 \ The general case }

\vskip10pt

  {\bf 1.1 Quantized universal enveloping algebras.} \, Let  $ \Bbbk $  be
a fixed field of zero characteristic, let  $ q $  be an indeterminate,
and let  $ \gerg $  be a semisimple Lie algebra over  $ \Bbbk \, $.
Let  $ \uqg $  be the quantum group \`a la Jimbo-Lusztig defined over
$ \kq \, $:  \, we define it after the conventions in [DP], or [DL],
or [Ga1] (for  $ \, \varphi = 0 \, $).  Actually, we can define a
quantum group like that for each lattice  $ M $  between the root
lattice  $ Q $  and the weight lattice  of  $ P $  of  $ \gerg \, $,
thus we shall write  $ \uqmg \, $.  Roughly,  $ \uqmg $  is the unital
$ \kq $--algebra  with generators  $ \, F_i \, $,  $ \Lambda_i^{\pm 1}
\, $,  $ E_i $  for  $ \, i=1, \dots, r =: \hbox{\it rank}\,(\gerg) \, $
and relations as in [DP], [Ga1], which depend on the Cartan datum
of  $ \gerg $  and on the choice of the lattice  $ M \, $;  \, in
particular, the  $ \Lambda_i $'s  are ``toral'' generators, roughly
$ q $--exponentials  of the elements of a  $ \Z $--basis  of  $ M $.
Here we only recall the relation
  $$  E_i F_j - F_j E_i \; = \; \delta_{i{}j} \, {{\, K_i - K_i^{-1} \,
\over \, q - q^{-1} \,}}   \qquad  \forall \;\; i, j = 1, \dots, r
\eqno (1.1)  $$
where  $ K_i $  is a  $ q $--analogue  of the coroot corresponding to
the  $ i $-th  node of the Dynkin diagram of  $ \gerg \, $  (in fact,
it is a suitable product of  $ \Lambda_k^{\pm 1} $'s).  Also, we
consider on  $ \uqmg $  the Hopf algebra structure given in [DP]
or [Ga1].
                                             \par
   The quantum Borel subalgebra  $ \uqmbp $  is simply the unital
$ \kq $--subalgebra  of  $ \uqmg $  generated by  $ \, \Lambda_1^{\pm 1}
\, $,  $ \dots \, $,  $ \Lambda_r^{\pm 1} \, $,  $ E_1 \, $,  $ \dots
\, $,  $ E_r \, $,  \, and  $ \uqmbm $  the one generated by  $ \, F_1
\, $,  $ \dots \, $,  $ F_r \, $,  $ \Lambda_1^{\pm 1} \, $,  $ \dots
\, $,  $ \Lambda_r^{\pm 1} \, $.  In fact, both of these are Hopf
$ \kq $--subalgebras  of  $ \uqmg \, $.  It follows that  $ \uqmg $
is generated by  $ \uqmbp $  and  $ \uqmbm \, $,  \, and every possible
commutation relation between these two subalgebras is a consequence of
(1.1) and the commutation relations between the  $ \Lambda_i^{\pm 1} $'s
and the  $ F_j $'s  or the  $ E_j $'s.  Finally, we call
$ \uqmt $  the unital  $ \kq $--subalgebra  of  $ \uqmg $  (and
of  $ U_q^{\scriptscriptstyle M}\!(\gerb_\pm) \, $)  generated
by all the  $ \Lambda_i $'s  ($ i = 1, \dots, n $),  which also
is a Hopf subalgebra.
                                               \par
   Mapping  $ \, F_i \mapsto E_i \, $,  $ \, \Lambda_i^{\pm 1}
\mapsto \Lambda_i^{\mp 1} \, $  and  $ \, E_i \mapsto F_i \, $
(for all  $ \, i = 1, \dots, n \, $)  \, uniquely defines an algebra
automorphism and coalgebra antiautomorphism of  $ \uqmg $,  that is a Hopf
algebra isomorphism  $ \; \varTheta \, \colon \, \uqmg \;{\buildrel \cong
\over {\lhook\joinrel\relbar\joinrel\relbar\joinrel\twoheadrightarrow}}\;
{\uqmg}^{\text{op}} \, $,  \; where hereafter given any Hopf algebra
$ H $  we denote by  $ H^{\text{op}} $  the same Hopf algebra as  $ H $
but for taking the opposite coproduct.  Restricting  $ \varTheta $
to quantum Borel subalgebras gives Hopf algebra isomorphisms  $ \;
U_q^{\scriptscriptstyle M}\!(\gerb_\pm) \cong U_q^{\scriptscriptstyle M}
(\gerb_\mp)^{\text{op}} \, $.

\vskip7pt

   {\bf 1.2 Quantum function algebras.} \, Let  $ M $  be a lattice between 
$ Q $  and  $ P $  as in \S 1.1, and define  $ \, M' := \big\{ \psi
\in P \,\big|\, \langle \psi, \mu \rangle \in \Z \, , \, \forall \, \mu
\in \Z \,\big\} \, $  where  $ \, \langle \ , \ \rangle \, $  is the 
$ \Q $--valued  scalar pro\-duct  on  $ P $  induced by scalar extension 
from the natural  $ \Z $--valued     
   \hbox{pairing between  $ Q $  and  $ P $.  Such}
  $ M' $  is again a lattice, said to be  {\sl dual\/}  to  $ M $. 
Conversely,  $ M $  is dual to  $ M' $,  i.e.~$ \, M = M'' \, $.
                                               \par
   We define quantum function algebras after Lusztig.  To start with,
letting  $ M $  and  $ M' $  be mutually dual lattices as above, we
define  $ \fqmpg $  as the unital  $ \kq $--algebra  of all matrix
coefficients of finite dimensional  $ \uqmg $--modules  which have
a basis of eigenvectors for all the  $ \Lambda_i $'s  ($ i = 1, \dots,
n $)  with eigenvalues powers of  $ q \, $.  Starting from  $ \uqmbp $
or  $ \uqmbm $  instead of  $ \uqmg $  the same recipe defines the
Borel quantum function algebras  $ \fqmpbp $  and  $ \fqmpbm $
respectively.
     \hbox{All these quantum function algebras are in fact
Hopf algebras too.}
                                               \par
   Finally, the Hopf algebra monomorphisms  $ \; j_\pm \, \colon \,
U_q^{\scriptscriptstyle M}\!(\gerb_\pm) \lhook\joinrel\longrightarrow
\uqmg \; $  induce Hopf algebra epimorphisms  $ \; \pi_\pm \, \colon
\, \fqmpg \relbar\joinrel\relbar\joinrel\twoheadrightarrow
F_q^{\scriptscriptstyle M'}\!\big[B_\pm\big] \; $.  See [DL]
and [Ga1] for details.

\vskip7pt

   {\bf 1.3 Isomorphisms between QUEA's and QFA's over Borel subgroups.}
\, Let  $ M $  and  $ M' $  be mutually dual lattices as in \S 1.2.
According to Tanisaki (cf.~[Ta]) there exist perfect (i.e.~non
degenerate) Hopf pairings  $ \; \uqmbp^{\text{op}} \! \otimes
\uqmpbm \longrightarrow \kq \, $,  $ \; \uqmbm^{\text{op}} \!
\otimes \uqmpbp \longrightarrow \kq \, $;  \; this implies  $ \;
\uqmbp^{\text{op}} \! \cong \fqmbm \; $  and  $ \; \uqmbm^{\text{op}}
\! \cong \fqmbp \, $.  Composing the latters with the isomorphisms
$ \; \uqmbp \cong \uqmbm^{\text{op}} \; $  and  $ \; \uqmbm \cong
\uqmbp^{\text{op}} \; $  in \S 1.1 it follows that  $ \; \uqmbp
\cong \fqmbp \; $  and  $ \; \uqmbm \cong \fqmbm \; $  as Hopf
$ \kq $--algebras.

\vskip7pt

   {\bf 1.4 Generation of  $ \uqmg $  by quantum function algebras.}
\, We said in \S 1.1 that  $ \uqmg $  is generated by  $ \uqmbm $  and
$ \uqmbp $,  whose mutual commutation is a consequence of (1.1).  In particular, we have a  $ \kq $--vector  space isomorphism  $ \; \uqmg
\, = \, \big( \uqmbm \otimes \uqmbp \big) \Big/ J \, $,  \; where  $ J $ 
is the two-sided ideal of  $ \, \uqmbm \otimes \uqmbp \, $   --- with
the standard tensor product structure ---   generated by  $ \, \big(
{\{ K_\mu \otimes 1 - 1 \otimes K_\mu \}}_{\mu \in M} \big) \, $,  \,
while the multiplication is a consequence of the internal commutation rules
of  $ U_q^{\scriptscriptstyle M}\!(\gerb_\pm) $  and by  (1.1).  Now, thanks
to the isomorphisms in \S 1.3, we describe  $ \uqmg $  as being generated
by  $ \fqmbm $  and  $ \fqmbp $,  with mutual commutation being a
consequence of the commutation formulas corresponding to (1.1) under
those isomorphisms.  So we have a  $ \kq $--vector  space isomorphism 
    $ \; \uqmg \cong $\allowbreak
  $ \, \big( \fqmbm \otimes \fqmbp \big) \Big/ I \, $,  \; where  $ I $
    \hbox{is the ideal of  $ \, \fqmbm \otimes \fqmbp \, $ 
corresponding to  $ J $,}   
while commutation rules are the internal ones of
$ F_q^{\scriptscriptstyle M}\big[B_\pm\big] $
and those corresponding to (1.1).

\vskip7pt

   {\bf 1.5 Relation with  $ L $--operators.} \, Tracking carefully
the construction of  $ \uqmg $  proposed in \S 1.4 above one realizes
that this is just an alternative way to introduce  $ \uqmg $  via
$ L $--operators  as made in [FRT].  Such a comparison is essentially
the meaning   --- or a possible interpretation ---   of the analysis
carried on in [Mo].  Moreover, the latter analysis also shows that
$ L $--operators  in [FRT] do correspond to suitable matrix coefficients
in  $ \fqmbm $  and  $ \fqmbp $  (embedded inside  $ \fqmg $);  such
matrix coefficients then correspond to quantum root vectors in
$ {\uqmbp}^{\text{op}} $  and  $ {\uqmbm}^{\text{op}} $  via the
isomorphisms  $ \, \fqmbm \cong {\uqmbp}^{\text{op}} \, $  and
$ \, \fqmbp \cong {\uqmbm}^{\text{op}} \, $  in \S 1.3, and finally
to quantum root vectors in  $ \uqmbm $  and  $ \uqmbp $  via the
isomorphisms  $ \, {\uqmbp}^{\text{op}} \cong \uqmbm \, $  and
$ \, {\uqmbm}^{\text{op}} \cong \uqmbp \, $  in \S 1.1.

\vskip7pt

   {\bf 1.6 Integer  $ \kqqm $--forms,  specializations, quantum
Frobenius morphisms.} \, In order to look at ``specializations of
a quantum group at special values of the parameter  $ q \, $'',
one needs the given quantum group to be defined over a subring of
$ \kq $  whose elements are regular, i.e.~have no poles, at such
special values.  As it is typical, we choose as ground ring the
Laurent polynomial ring  $ \kqqm $.  Then instead of  $ \uqmg $
we must consider integer forms of  $ \uqmg $  over  $ \kqqm $,
i.e.~Hopf  $ \kqqm $--subalgebras  of  $ \uqmg $
which give back all of  $ \uqmg $  by scalar extension from  $ \kqqm $
to  $ \kq $:  if  $ \overline{U}_{\! q}^{\scriptscriptstyle \, M}\!
(\gerg) $  is such a  $ \kqqm $--form,  its  {\sl specialization\/}
at  $ \, q = c \in \Bbbk \, $  is the quotient Hopf  $ \Bbbk $--algebra
$ \; \overline{U}_{\! c}^{\scriptscriptstyle \, M}\!(\gerg) :=
\overline{U}_{\! q}^{\scriptscriptstyle \, M}\!(\gerg) \! \Big/ \!
(q - c) \, \overline{U}_{\! q}^{\scriptscriptstyle \, M}\!(\gerg) \; $.
                                                      \par
   There are essentially two main types of  $ \kqqm $--integer  forms:
one is  $ \uhat_q^{\scriptscriptstyle M}\!(\gerg) $  (the quantum
analogue of Kostant's  $ \Z $--integer  form of  $ \gerg \, $)
introduced by Lusztig in [Lu1], generated by  $ q $--binomial
coefficients and  $ q $--divided  powers); the second one is
$ \utilde_q^{\scriptscriptstyle M}\!(\gerg) \, $,  introduced
by De Concini and Procesi in [DP], generated by rescaled
quantum root vectors; see [Ga1] for details.  When  $ q $
is specialized to any value in  $ \Bbbk $  {\sl which is not
a root of 1},  the choice of either of these two integer forms
is irrelevant, because the corresponding specialized Hopf
$ \Bbbk $--algebras  are mutually isomorphic.  If instead  $ q $
is specialized to  $ \, \varepsilon \in \Bbbk \, $  {\sl which is
a root of 1},  then the specialized algebra changes according to
the choice of integer form.
                                                      \par
   Indeed, the behavior of  $ \uhat_q^{\scriptscriptstyle M}\!
(\gerg) $  and  $ \utilde_q^{\scriptscriptstyle M}\!(\gerg) $
w.r.t.~specializations at roots of 1 is pretty different,
even opposite.  In particular, one has semiclassical limits
$ \; \uhat_1^{\scriptscriptstyle M}\!(\gerg) \cong U(\gerg) \, $,
\, the universal enveloping algebra of  $ \gerg \, $,  \, and
$ \; \utilde_1^{\scriptscriptstyle M}\!(\gerg) \cong F \big[
G^*_{\scriptscriptstyle M} \big] \, $,  \, the regular function
algebra of   $ G^*_{\scriptscriptstyle M} $,  where
$ G^*_{\scriptscriptstyle M} $  is a connected Poisson
algebraic group with fundamental group isomorphic to  $ \,
P \big/ M \, $  and dual to  $ \gerg \, $,  \, the latter
endowed with a structure of Lie bialgebra, inherited from
$ \uhat_q^{\scriptscriptstyle M}\!(\gerg) \, $.  Moreover,
specializations of an integer form of either type at a root of 1,
say  $ \, \varepsilon \in \Bbbk \, $,  are linked to its semiclassical
limit by the so-called  {\sl quantum Frobenius morphisms}
  $$  \uhat_\varepsilon^{\scriptscriptstyle M}\!(\gerg) \,
\relbar\joinrel\relbar\joinrel\relbar\joinrel\twoheadrightarrow
\, \uhat_1^{\scriptscriptstyle M}\!(\gerg) \, \cong \, U(\gerg)
\quad ,  \qquad
   F\big[G^*_{\scriptscriptstyle M}\big] \cong \,
\utilde_1^{\scriptscriptstyle M}\!(\gerg) \,
\lhook\joinrel\relbar\joinrel\relbar\joinrel\relbar\joinrel\rightarrow
\, \utilde_\varepsilon^{\scriptscriptstyle M}\!(\gerg) \quad .
\eqno (1.2)  $$
   \indent   Such a situation occurs exactly the same   --- {\sl mutatis
mutandis}  ---   for the quantum Borel subalgebras  $ \uqmbm $  and
$ \uqmbp \, $.  In short, one has two types of  $ \kqqm $--integer
forms  $ \uhat_q^{\scriptscriptstyle M}\!(\gerb_\pm) $  and
$ \utilde_q^{\scriptscriptstyle M}\!(\gerb_\pm) $,  and quantum
Frobenius morphisms
  $$  \uhat_\varepsilon^{\scriptscriptstyle M}\!(\gerb_\pm)
\, \relbar\joinrel\relbar\joinrel\relbar\joinrel\twoheadrightarrow \,
\uhat_1^{\scriptscriptstyle M}\!(\gerb_\pm) \, \cong \, U(\gerb_\pm)
\quad ,  \qquad
   F\big[B_\pm^*\big] \cong \,
\utilde_1^{\scriptscriptstyle M}\!(\gerb_\pm) \,
\lhook\joinrel\relbar\joinrel\relbar\joinrel\relbar\joinrel\rightarrow
\, \utilde_\varepsilon^{\scriptscriptstyle M}\!(\gerb_\pm) \quad .
\eqno (1.3)  $$
By construction,  $ \uhat_q^{\scriptscriptstyle M}\!(\gerg) $  is
generated by  $ \uhat_q^{\scriptscriptstyle M}\!(\gerb_+) $  and
$ \uhat_q^{\scriptscriptstyle M}\!(\gerb_-) $,  and similarly
$ \utilde_q^{\scriptscriptstyle M}\!(\gerg) $  is generated
by  $ \utilde_q^{\scriptscriptstyle M}\!(\gerb_+) $  and
$ \utilde_q^{\scriptscriptstyle M}\!(\gerb_-) $.  It follows
that the morphisms in (1.3) can also be obtained from (1.2) by
restriction to quantum Borel subalgebras; conversely, the quantum
Frobenius morphisms in (1.2) are uniquely determined   --- and
described ---   by those in (1.3).
                                            \par
   By duality, the like happens also for quantum function algebras:
in particular, there exist two  $ \kqqm $--integer  forms
$ \fhat_q^{\scriptscriptstyle M}\!\big[G\big] $  and
$ \ftilde_q^{\scriptscriptstyle M}\!\big[G\big] $  of  $ \fqmg $,
which are dual respectively to  $ \uhat_q^{\scriptscriptstyle M}
\!(\gerg) $  and  $ \utilde_q^{\scriptscriptstyle M}\!(\gerg) $
in Hopf theoretical sense, for which the dual of (1.2) holds, namely
  $$  F\big[G\big] \cong \,
\fhat_1^{\scriptscriptstyle M}\!\big[G\big] \,
\lhook\joinrel\relbar\joinrel\relbar\joinrel\relbar\joinrel\rightarrow
\, \fhat_\varepsilon^{\scriptscriptstyle M}\!\big[G\big]
\quad ,  \qquad
   \ftilde_\varepsilon^{\scriptscriptstyle M}\!\big[G\big] \,
\relbar\joinrel\relbar\joinrel\relbar\joinrel\twoheadrightarrow
\, \ftilde_1^{\scriptscriptstyle M}\!\big[G\big] \cong \,
U(\gerg^*)  \quad .   \eqno (1.4)  $$
Similarly, the dual of (1.3) holds for quantum function algebras
of Borel subgroups, namely
  $$  F\big[B_\pm\big] \cong \,
\fhat_1^{\scriptscriptstyle M}\!\big[B_\pm\big] \,
\lhook\joinrel\relbar\joinrel\relbar\joinrel\relbar\joinrel\rightarrow
\, \fhat_\varepsilon^{\scriptscriptstyle M}\!\big[B_\pm\big]
\quad ,  \qquad
   \ftilde_\varepsilon^{\scriptscriptstyle M}\!\big[B_\pm\big]
\, \relbar\joinrel\relbar\joinrel\relbar\joinrel\twoheadrightarrow
\, \ftilde_1^{\scriptscriptstyle M}\!\big[B_\pm\big] \cong
\, U(\gerb_\pm^*)  \quad ,   \eqno (1.5)  $$
which follow from (1.4) via the maps  $ \, \fqmg \, {\buildrel
{\pi_\pm} \over {\relbar\joinrel\relbar\joinrel\twoheadrightarrow}}
\, F_q^{\scriptscriptstyle M}\!\big[B_\pm\big] \, $  in \S 1.2.
See [Ga1] for details.
                                              \par
   The point we want to stress now is the relation between the
isomorphisms of Hopf  $ \kq $--algebras  $ \; \uqmbp \cong \fqmbp
\; $  and  $ \; \uqmbm \cong \fqmbm \; $  in \S 1.3 and the
$ \kqqm $--integer  forms on both sides.  The key fact is that
the previous isomorphisms restrict to isomorphisms of Hopf
$ \kqqm $--algebras  $ \; \uhat_q^{\scriptscriptstyle M}\!
(\gerb_\pm) \cong \ftilde_q^{\scriptscriptstyle M}\!\big[
B_\pm \big] \; $  and  $ \; \utilde_q^{\scriptscriptstyle M}\!
(\gerb_\pm) \cong \fhat_q^{\scriptscriptstyle M}\!\big[B_\pm\big] \; $.
Therefore, looking at  $ \uqmg $  as generated by  $ \fqmbm $  and
$ \fqmbp $  as explained in \S 1.4 one argues that the  {\sl first},
resp.~the  {\sl second},  quantum Frobenius morphisms in (1.2) are
uniquely determined (and described) by the   {\sl second\/}  ones,
resp.~the  {\sl first\/}  ones, in (1.5).

\vskip1,3truecm

\centerline {\bf \S \; 2 \ The case of  $ \gergln \, $ }

\vskip10pt

{\bf 2.1  $ q $--matrices.}  Let  $ \big\{ t_{i{}j} \,\big|\, i, j
= 1, \ldots, n \big\} $  be a set of elements in any  $ \kq $--algebra
$ A \, $,  ideally displayed inside an  $ (n \times n) $--matrix
they are the entries of.  We'll say that  $ \, T := {\big( t_{i{}j}
\big)}_{i,j=1,\dots,n} \, $  {\sl is a  $ q $--matrix\/}  if
the  $ t_{i{}j} $'s  enjoy the following relations
  $$  \hbox{ $ \eqalign {
   {} \hskip11pt   t_{ij} \, t_{ik} = q \, t_{ik} \, t_{ij} \; ,
\quad \quad  t_{ik} \, t_{hk}  &  = q \, t_{hk} \, t_{ik}
\hskip56pt  \forall  \quad  j<k \, , \, i<h \, ,  \cr
   {} \hskip11pt   t_{il} \, t_{jk} = t_{jk} \, t_{il} \; ,
\qquad  t_{ik} \, t_{jl} - \, t_{jl} \, t_{ik}  &
= \left( q - q^{-1} \right) \, t_{il} \, t_{jk}   \hskip20pt
\forall  \quad  i<j \, , \, k<l \, .  \cr } $ }   \eqno  $$
in the algebra  $ A $.  {\sl In this case, the so-called  {\it
``quantum determinant''},  defined as
  $$  \hbox{\sl det}_q \left(\! {\big(t_{k,\ell}\big)}_{k, \ell = 1,
\dots, n} \right) \, := {\textstyle \sum_{\sigma \in \Cal{S}_n}}
{(-q)}^{l(\sigma)} t_{1,\sigma(1)} \, t_{2,\sigma(2)} \cdots
t_{n,\sigma(n)}  $$
commutes with all the  $ t_{i,j} $'s.}  If in addition  $ A $  is
a  $ \kq $--bialgebra, we shall also require that   
  $$  \Delta (t_{ij}) \, = \, {\textstyle \sum\limits_{k=1}^n}
\, t_{ik} \otimes t_{kj} \; ,   \qquad   
      \epsilon(t_{ij}) \, = \, \delta_{ij}  
\eqno   \forall \;\; i, j = 1, \dots, n \; . \qquad  $$  
In this case, the quantum determinant is group-like,
that is  $ \, \Delta(\text{\sl det}_q) = \text{\sl det}_q \otimes
\text{\sl det}_q \, $  and  $ \, \epsilon(\text{\sl det}_q) = 1
\, $.  Finally, if  $ A $  is a Hopf algebra we call  {\sl  Hopf 
$ q $--matrix\/}  any  $ q $--matrix  like above whose entries
are such that  $ \hbox{\sl det}_q $  is invertible in  $ A \, $; 
\, then  $ \, S \big( \text{\sl det}_q^{\, \pm 1} \big) =
\text{\sl det}_q^{\, \mp 1} \, $.
                                        \par
   For later use we also recall the following compact notation.
Let  $ \, T_1 := T \otimes I \, $,  $ \, T_2 := I \otimes T \in
A \otimes {\text{\sl Mat}_n\!\big(\kq\big)}^{\otimes 2} \cong A
\otimes \text{\sl Mat}_{n^2}\!\big(\kq\big) \, $,  \, where  $ I $
is the identity matrix, and  $ \, T := {\big( t_{i{}j} \big)}_{i,j=1,
\dots,n} \, $  is thought of as an element of  $ \, \text{\sl Mat}_n
\!\big(A\big) \cong A \otimes \text{\sl Mat}_n\!\big(\kq\big) \, $;
\, consider
  $$  R \, := \, \textstyle{\sum_{i,j=1}^n} \, q^{\delta_{i{}j}}
\, e_{i{}i} \otimes e_{j{}j} \, + \, \big( q - q^{-1} \big) \,
\textstyle{\sum_{1 \leq i < j \leq n}} \, e_{i{}j} \otimes
e_{j{}i} \; \in \; \text{\sl Mat}_{n^2}\!\big(\kq\big)  $$
where  $ \, e_{i{}j} := \big( \delta_{i{}h} \, \delta_{j{}k}
\big)_{h,k=1}^n \, $  is the  $ (i,j) $--th  elementary matrix.
Then  $ T $  is a  $ q $--matrix  if and only if the identity
$ \; R \, T_2 \, T_1 = T_1 \, T_2 \, R \; $  holds true in
$ \, A \otimes \text{\sl Mat}_{n^2}\!\big(\kq\big) \, $;  \,
in detail, for the matrix entry in position  $ \big( (i,j),
(kl) \big) $  this reads  $ \; \sum\limits_{m,p=1}^n R_{ij,mp}
\; t_{pk} \, t_{ml} \, = \sum\limits_{m,p=1}^n  t_{im} \, t_{jp}
\; R_{mp,kl} \; $.  In the bialgebra case  $ T $  is a  $ q $--matrix 
if in addition  $ \, \Delta(T) = T \,\dot\otimes\, T \, $,  $ \,
\epsilon(T) = I \, $,  \, and in the Hopf algebra case also  $ \,
T \, S(T) = I = S(T) \, T \, $,  \, i.e.~$ \, S(T) = T^{-1} \, $; 
\, see [FRT] and [No] for notations   --- we use assumptions and
normalizations of the latter ---   and further details.

\vskip7pt

  {\bf 2.2 Presentation of  $ \fqpg $,  $ \fqpbm $  and  $ \fqpbp $
for  $ \, G = GL_n \, $.} \, Let's look at  $ \, G = GL_n \, $.
After [APW], Appendix, we know that  $ \fqpgl $  has the following
presentation: it is the unital associative  $ \kq $--algebra  with
generators the elements of  $ \; \big\{ t_{ij} \,\big|\, i, j = 1,
\ldots, n \big\} \bigcup \big\{ \text{\sl det}_q^{\, -1} \big\} \; $
and relations encoded by the requirement that  $ \, {\big( t_{i,j}
\big)}_{i,j=1,\dots,n} \, $  be a  $ q $--matrix;  in particular,
$ \text{\sl det}_q^{\, \pm 1} $  belongs to the centre of  $ \fqpgl $.
Moreover,  $ \fqpgl $  has the unique Hopf algebra structure such
that  $ \, {\big( t_{i,j} \big)}_{i, j = 1, \dots, n} \, $
be a Hopf $ q $--matrix.
                                                \par
   Similarly,  $ \fqpbm $  and  $ \fqpbp $  are defined in the same
way  {\sl but\/}  with the additional relations  $ \; t_{i,j} = 0 \,
(i, j = 1, \dots, n; i > j) \; $  for  $ \fqpbm $  and  $ \; t_{i,j}
= 0 \, (i, j = 1, \dots, n; i < j) \; $  for  $ \fqpbp \, $.  Otherwise,
we can say that  $ \fqpbm $,  respectively  $ \fqpbp $,  is generated
by the entries of the  $ q $--matrix
  $$  \pmatrix
  t_{1,1}  &  0  &  \cdots  &  0  &  0   \\
  t_{2,1}  &  t_{2,2}  &  \cdots  &  0  &  0   \\
  \vdots  &  \vdots  &  \vdots  &  \vdots  &  \vdots   \\
  t_{n-1,1}  &  t_{n-1,2}  &  \cdots  &  t_{n-1,n-1}
&  0   \\
  t_{n,1}  &  t_{n,2}  &  \cdots  &  t_{n,n-1}
&  t_{n,n}   \\
   \endpmatrix
\, ,  \, \text{resp.} \,
      \pmatrix
  t_{1,1}  &  t_{1,2}  &  \cdots  &  t_{1,n-1}
&  t_{1,n}   \\
  0  &  t_{2,2}  &  \cdots  &  t_{2,n-1}  &  t_{2,n}   \\
  \vdots  &  \vdots  &  \vdots  &  \vdots  &  \vdots   \\
  0  &  0  &  \cdots  &  t_{n-1,n-1}  &  t_{n-1,n}   \\
  0  &  0  &  \cdots  &  0  &  t_{n,n}   \\
   \endpmatrix  $$
and by the additional element  $ \, \big( t_{1,1} \, t_{2,2}
\cdots t_{n,n} \big)^{-1} \, $.  Moreover, both  $ \fqpbm $  and
$ \fqpbp $  are Hopf algebras, the Hopf structure being given by
the assumption that their generating matrices be  {\sl Hopf}
$ q $--matrices.  See also [PW] for all these definitions.
                                                  \par
   By the very definitions, the Hopf algebra epimorphisms  $ \;
\pi_+ \, \colon \, \fqpgl \longtwoheadrightarrow \fqpbp \; $  and
$ \; \pi_- \, \colon \, \fqpgl \longtwoheadrightarrow \fqpbm \; $
mentioned in \S 1.2 are given by  $ \; \pi_+ : \, t_{ij} \mapsto
t_{ij} \; (i \leq j) \, $,  $ \; t_{ij} \mapsto 0 \; (i > j)
\; $  and  $ \; \pi_- : \, t_{ij} \mapsto t_{ij} \; (i \geq j)
\, , \; t_{ij} \mapsto 0 \; (i < j) \; $  respectively.

\vskip7pt

  {\bf 2.3. The quantum algebras  $ \uqmg $,  $ \uqmbm $  and
$ \uqmbp $  for  $ \, \gerg = \gergln \, $,  $ \,
M \in \{ P, Q \} \, $.} \  We recall (cf.~for instance [GL]) the
definition of the quantized universal enveloping algebra  $ \uqpgl $:
it is the associative algebra with 1 over  $ \kq $  with generators
  $$  F_1 \, , \, F_2 \, , \, \dots \, , \, F_{n-1} \, , \, G_1^{\pm 1}
\, , \, G_2^{\pm 1} \, , \, \dots \, , \, G_{n-1}^{\pm 1} \, , \,
G_n^{\pm 1} \, , \, E_1 \, , \, E_2 \, , \, \dots \, , \, E_{n-1}  $$
and relations
  $$  \displaylines {
   \hfill   G_i G_i^{-1} = \, 1 = \, G_i^{-1} G_i \; ,  \qquad
G_i^{\pm 1} G_j^{\pm 1} = \, G_j^{\pm 1} G_i^{\pm 1}  \qquad
\hfill  \forall \, i, j \qquad  \cr
   \hfill   G_i F_j G_i^{-1} = \, q^{\delta_{i,j+1} - \delta_{i,j}} F_j
\; ,  \qquad  G_i E_j G_i^{-1} = \, q^{\delta_{i,j} - \delta_{i,j+1}}
E_j  \qquad  \hfill  \forall \, i, j \qquad  \cr
   \hfill   E_i F_j - F_j E_i \; = \; \delta_{i,j} \, {{\; G_i G_{i+1}^{-1}
- G_i^{-1} G_{i+1} \;} \over {\; q - \qm \;}}  \qquad \qquad   \hfill
\forall \, i, j  \qquad  \cr
 }  $$   
  $$  \displaylines {
   \hfill   \quad  E_i E_j \, = \, E_j E_i \; ,  \qquad  F_i F_j \,
= \, F_j F_i   \hfill  \forall \; i, j : \vert i - j \vert > 1
\, \phantom{.} \;  \cr
   \hfill   E_i^2 E_j - [2]_q E_i E_j E_i + E_j E_i^2 \, = \, 0 \; ,
\quad  F_i^2 F_j - [2]_q F_i F_j F_i + F_j F_i^2 \, = \, 0
\hfill  \forall \; i, j : \vert i - j \vert = 1  \cr }  $$
with  $ \, [2]_q := q + q^{-1} \, $.  Moreover,  $ \uqpgl $  has
a Hopf algebra structure given by
  $$  \displaylines {
   \hfill   \Delta \left( F_i \right) = F_i \otimes G_i^{-1} G_{i+1}
+ 1 \otimes F_i \, ,  \; \qquad  S \left( F_i \right) = - F_i G_i
G_{i+1}^{-1} \, , \; \qquad  \epsilon \left( F_i \right) = 0
\hfill  \forall \; i \, \phantom{.} \;  \cr
   \hfill   \qquad  \Delta \left( G_i^{\pm 1} \right) = G_i^{\pm 1}
\otimes G_i^{\pm 1} \, ,  \qquad \qquad \qquad  S \left( G_i^{\pm 1}
\right) = G_i^{\mp 1} \, ,  \quad \qquad  \epsilon \left( G_i^{\pm 1}
\right) = 1   \hfill  \forall \; i \, \phantom{.} \;  \cr
   \hfill   \Delta \left( E_i \right) = E_i \otimes 1 + G_i G_{i+1}^{-1}
\otimes E_i \, ,  \qquad  S \left( E_i \right) = - G_i^{-1} G_{i+1} E_i
\, ,  \qquad  \epsilon \left( E_i \right) = 0   \hfill  \forall \; i
\, . \;  \cr }  $$
   \indent   The algebra  $ \uqqgl $   --- defined as in [Ga1], \S 3 ---
can be realized as a Hopf subalgebra.  Namely, define  $ \; L_i := G_1
G_2 \cdots G_i \, $,  $ \, K_j := G_j G_{j+1}^{-1} \, $  for all  $ \,
i = 1, \dots, n \, $,  $ \, j = 1, \dots, n-1 \, $.  Then  $ \uqqgl $
is the  $ \kq $--subalgebra  of  $ \uqpgl $  generated by  $ \, \big\{
F_1, \dots, F_{n-1}, K_1^{\pm 1}, \dots,
 \allowbreak
K_{n-1}^{\pm 1}, L_n^{\pm 1},
E_1, \dots, E_{n-1} \big\} \, $.
   The quantum Borel subalgebra  $ \uqpbp $,  resp.~$ \uqpbm $,  is the
subalgebra of  $ \uqpgl $  generated by  $ \, \big\{ G_1^{\pm 1}, \dots,
G_n^{\pm 1} \big\} \cup \big\{ E_1, \dots, E_{n-1} \big\} \, $,  resp.~by
$ \, \big\{ G_1^{\pm 1}, \dots,
 \allowbreak
G_n^{\pm 1} \big\} \cup \big\{ F_1,
\dots, F_{n-1} \big\} \, $.  Similar definitions hold for
$ U_q^{\scriptscriptstyle Q}\!(\gerb_\pm) $,  but with the set
$ \big\{ K_1^{\pm 1}, \dots,
 \allowbreak
K_{n-1}^{\pm 1}, L_n^{\pm 1} \big\} \, $
instead of  $ \big\{ G_1^{\pm 1}, \dots, G_n^{\pm 1} \big\} \, $.
All these are in fact  {\sl Hopf}  subalgebras.

\vskip7pt

  {\bf 2.4. The Hopf isomorphisms  $ \, \zeta_- \colon \uqpbm \cong
\fqpbm \, $,  $ \, \zeta_+ \colon \uqpbp \cong \fqpbp \, $.} \, The
Hopf algebra isomorphisms of \S 1.3 are given explicitly by  ($ \,
i=1, \dots, n; j=1, \dots, n-1 \, $)
  $$  \displaylines{
   \zeta_- : \uqpbm \,{\buildrel \cong \over \loongrightarrow}\,
\fqpbm \; ,  \qquad
  G_i^{\pm 1} \mapsto t_{i,i}^{\; \mp 1} \; ,  \quad
F_j \mapsto + \big( q - q^{-1} \big)^{-1} t_{j+1,j+1}^{\;-1}
\, t_{j+1,j}  \cr
   \zeta_+ : \uqpbp \,{\buildrel \cong \over \loongrightarrow}\,
\fqpbp \; ,  \qquad
  G_i^{\pm 1} \mapsto t_{i,i}^{\; \pm 1} \; ,  \quad
E_j \mapsto - \big( q - q^{-1} \big)^{-1} t_{j,j+1}
\, t_{j+1,j+1}^{\;-1}  \cr }  $$
and their inverse are uniquely determined by
  $$  \displaylines{
   \zeta_-^{\,-1} : \fqpbm \,{\buildrel \cong \over \loongrightarrow}\,
\uqpbm \; ,  \qquad
  t_{i,i}^{\; \pm 1} \mapsto G_i^{\mp 1} \; ,  \quad
t_{j+1,j} \mapsto + \big( q - q^{-1} \big) \, G_{j+1}^{\;-1} \, F_j
\; \phantom{.}  \cr
   \zeta_+^{\,-1} : \fqpbp \,{\buildrel \cong \over \loongrightarrow}\,
\uqpbp \; ,  \qquad
  t_{i,i}^{\; \pm 1} \mapsto G_i^{\pm 1} \; ,  \quad
t_{j,j+1} \mapsto - \big( q - q^{-1} \big) \, E_j \,
G_{j+1}^{\; +1} \; .  \cr }  $$
A straightforward computation shows that all the above are
isomorphisms as claimed.

\vskip7pt

\proclaim{Theorem 2.5}  {\sl (``short'' FRT-like presentation of
$ \uqpgl \, $)}
                                         \hfill\break
   \indent   $ \, \uqpgl \, $  is the unital associative
$ \kq $--algebra  with generators the elements of the set  $ \;
\displaystyle{ \big\{\beta_{i,j}\big\}_{1 \leq i \leq j \leq n} \,
\bigcup \, \big\{\gamma_{j,i}\big\}_{1 \leq i \leq j \leq n}} \; $
and relations
  $$  \hbox{ $  \eqalign{
  \beta_{i,i+1} \, \gamma_{j+1,j} \, - \, \gamma_{j+1,j} \, \beta_{i,i+1}
\;  = \;  &  \big( \delta_{i,j+1} \, \big( 1 - q^{-1} \big) \, + \,
\delta_{i,j-1} \, (1 - q) \big) \, \beta_{i,i+1} \, \gamma_{j+1,j}
\; -  \cr
          &  \qquad \qquad  - \, \delta_{i{}j} \, \big( q - q^{-1} \big)
\big( \alpha_i \, \alpha_{i+1}^{-1} \, - \, \alpha_i^{-1} \,
\alpha_{i+1} \big)  \cr } $ }   \eqno (2.1)  $$
 \vskip-9pt
  $$  \beta_{k,k} \, \gamma_{k,k} = 1   \eqno (2.2)  $$
(for all  $ \, i, j = 1, \dots, n-1 \, $,  $ \, k = 1, \dots,
n \, $)  plus the relations encoded in the requirement that the
triangular matrices  $ \, B := \big( \beta_{i{}j} \big)_{i,j=1}^n
\, $  and  $ \, \varGamma := \big( \gamma_{i{}j} \big)_{i,j=1}^n
\, $  be  $ q $--matrices.  Moreover, this algebra has the unique
Hopf algebra structure such that these are Hopf  $ q $--matrices.
\endproclaim

\demo{Proof} This follows directly from \S 1.4 and the
isomorphisms in \S 2.4.  Indeed, in the given presentation the
$ \beta_{h,k} $'s  generate a copy of  $ \fqpbp $,  with  $ \,
\beta_{h,k} \cong t_{h,k} \, $,  isomorphic to  $ \uqpbp $
via \S 2.4; similarly, the  $ \gamma_{r,s} $'s  generate a copy
of  $ \fqpbm $,  with  $ \, \gamma_{r,s} \cong t_{r,s} \, $,
isomorphic to  $ \uqpbm $.  The additional set of ``mixed''
relations (2.1) involving simultaneously the  $ \beta_{i,i+1} $'s
and the  $ \gamma_{j+1,j} $'s  then correspond to the set of
relations (1.1)   --- or to the third line of the set of relations
in \S 2.3 ---   via the isomorphisms  $ \zeta_\pm $  of \S 2.4;
indeed, these isomorphisms give  $ \, \beta_{i,i+1} \cong - \big(
q - \qm \big) E_i \, G_{i+1}^{\; +1} \, $,  $ \, \beta_{k,k} \cong
G_k \, $,  and  $ \, \gamma_{j+1,j} \cong + \big( q - \qm \big)
G_{j+1}^{\; -1} \, F_j \, $,  $ \, \gamma_{k,k} \cong G_k^{\,-1}
\, $,  \, whence computing  $ \; - \big( q - \qm \big)^{\! 2} \,
\big[ E_i \, G_{i+1}^{\; +1} \, , G_{j+1}^{\; -1} \, F_j \big] \; $
in  $ \uqpgl \, $  we get formula (2.1).  As to the Hopf structure,
it is determined by that of the Hopf subalgebras  $ \uqpbp $  and
$ \uqpbm $:  thus the claim follows from the previous discussion.
\qed
\enddemo

\vskip7pt

   {\bf 2.6 Remark:}  note that any other commutation relation
between a generator  $ \beta_{h,k} $  ($ h < k $)  and a generator
$ \gamma_{r,s} $  ($ r > s $)  can be deduced from the ones between
the  $ \beta_{i,i+1} $'s  and the  $ \gamma_{j+1,j} $'s  using
repeatedly the relations
  $$  \beta_{i,j} \; = \; \big( q - q^{-1} \big)^{-1} \, \big(
\, \beta_{i,k} \, \beta_{k,j} \, - \, \beta_{k,j} \, \beta_{i,k}
\big) \; \beta_{k,k}^{\;-1}  \qquad \quad  (\, \forall \;\, i <
k < j \,)  $$
which spring out of the relations  $ \; \beta_{i,k} \, \beta_{k,j}
\, - \, \beta_{k,j} \, \beta_{i,k} \, = \, \big( q - q^{-1} \big)
\, \beta_{k,k} \, \beta_{i,j} \; $  for the  $ q $--matrix  $ B \, $,
\, and the relations
  $$  \gamma_{j,i} \; = \; \big( q - q^{-1} \big)^{-1} \, \big( \,
\gamma_{k,i} \, \gamma_{j,k} \, - \, \gamma_{j,k} \, \gamma_{k,i}
\big) \; \gamma_{k,k}^{\;+1}  \qquad \quad  (\, \forall \;\, j >
k > i \,)  $$
which arise from the relations  $ \; \gamma_{k,i} \, \gamma_{j,k} \,
- \, \gamma_{j,k} \, \gamma_{k,i} \, = \, \big( q - q^{-1} \big) \,
\gamma_{k,k} \, \gamma_{j,i} \; $
     \hbox{for the  $ q $--matrix  $ \varGamma \, $.}

\vskip7pt

  {\bf 2.7 Quantum root vectors and  $ L $--operators.} \,
In this subsection we describe the generators of  $ \uqpgl $
considered in Theorem 2.5 in terms of generators of the
Faddeev-Reshetikhin-Takhtajan (FRT in short) presentation,
the so-called  $ L $--operators,  ---   in [FRT].
                                           \par
   Our comparison ``factors through'' that with quantum root
vectors built upon the Jimbo-Lusztig generators given in \S 2.3.
For any  $ \, x $,  $ y $,  $ a $,  let  $ \, [w,y]_a := x \, y
- a \, y \, x \, $.  Define
  $$  \matrix
   E^\pm_{i,i+1} := E_i \; ,  &  \quad  E^\pm_{i,j} :=
\big[\, E^\pm_{i,k} \, , E^\pm_{k,j} \,\big]_{q^{\pm 1}}  &
\qquad  \forall \; i < k < j  \\
   F^\pm_{i+1,i} := F_i \; ,  &  \quad  F^\pm_{j,i} :=
\big[\, F^\pm_{j,k} \, , F^\pm_{k,i} \,\big]_{q^{\mp 1}}  &
\qquad  \forall \; j > k > i  \endmatrix  $$
as in [Ji]: all these are quantum root vectors, in that in the
semiclassical limit at  $ \, q = 1 \, $  they specialize to root
vectors for  $ \gergl_n \, $,  \, namely the elementary matrices
$ e_{i{}j} $  with  $ \, i \not= j \, $.  As a matter of notation,
     \hbox{set also  $ \, \edot_{i,j}^\pm := \big( q - q^{-1} \big) \,
E_{i,j}^\pm \, $  and  $ \, \fdot_{j,i}^\pm := \big( q - q^{-1} \big)
\, F_{j,i}^\pm \, $  for all  $ \, i < j \, $.}
                                           \par
   For the  $ L $--operators,  introduced in [FRT], we recall from
[No], \S 1.2, the formulas
  $$  \hbox{ $ \matrix
   L^+_{i{}i} := G_i^{\,+1} \; ,  &  \quad  L^+_{i{}j} :=
+ G_i^{\,+1} \, \fdot_{j,i}^+ \; ,  &  \quad  L^+_{j,i} := 0
&  \qquad  \forall \;\; i < j  \\
   L^-_{i{}i} := G_i^{\,-1} \; ,  &  \quad  L^-_{j{}i} :=
- \edot_{i,j}^+ \, G_i^{\,-1} \; ,  &  \quad  L^-_{i,j} := 0
&  \qquad  \forall \;\; i < j  \endmatrix $ }   \eqno (2.3)  $$
to define them; setting  $ \; L^+ := \big( L^+_{i{}j} \big)_{i,j=1}^n
\; $  and  $ \; L^- := \big( L^-_{i{}j} \big)_{i,j=1}^n \, $,  \; the
relations
  $$  R \, L^+_1 \, L^+_2 \, = \, L^+_2 \, L^+_1 \, R \; ,  \qquad
R \, L^-_1 \, L^-_2 \, = \, L^-_2 \, L^-_1 \, R \; ,  \qquad
R \, L^+_1 \, L^-_2 \, = \, L^-_2 \, L^+_1 \, R   \eqno (2.4)  $$
express in compact form their mutual commutation properties (with
notation as in \S 2.1).  Indeed, the FRT presentation amounts exactly
to claim that  $ \uqpgl $  is the unital associative  $ \kq $--algebra
with generators  $ L^\pm_{i,j} $  (for all  $ \, i, j = 1, \dots, n \, $)
and relations (2.4) and
  $$  L^+_{k,k} \, L^-_{k,k} \; = \; 1 \; = \; L^-_{k,k} \, L^+_{k,k}
\qquad  \forall \; k = 1, \dots, n   \eqno (2.5)  $$
and it has the unique Hopf algebra structure such that
  $$  \Delta(L^\varepsilon) \, = \, L^\varepsilon \,\dot\otimes\,
L^\varepsilon \; ,  \qquad  \epsilon(L^\varepsilon) \, = \, I \; ,
\qquad  S(L^\varepsilon) = \big(L^\varepsilon\big)^{-1}  \qquad
\quad  \forall \;\; \varepsilon \in \big\{ + \, , - \big\}
\eqno (2.6)  $$
where  $ L^+ $  and  $ L^- $  are the upper or lower triangular matrices
whose non-zero entries are the  $ L^+_{i,j} $'s  and the  $ L^-_{j,i} $'s
respectively,  $ I $  is the  $ (n \times n) $--identity  matrix and we
use standard compact notation as in [FRT] or [CP].
                                           \par
   Now, using the identifications  $ \zeta_+^{\,\pm 1} $  we get
identities
  $$  \beta_{i,i} \, = \, G_i^{\,+1} \; ,  \qquad
\beta_{i,j} \, = \, + {(-q)}^{j-i} \, G_j^{\,+1} \, \edot_{i,j}^-
\eqno \forall \;\; i < j \; . \quad \qquad (2.7)  $$
Indeed, the identities  $ \, \beta_{i{}i} = G_i^{\,+1} \, $  and  $ \,
\beta_{i,j} = - q \, G_j^{\,+1} \, \edot_{i,j}^- = - \, \edot_{i,j}^-
\, G_j^{\,+1} \, $  for  $ \, j = i + 1 \, $  come out directly from
the description of  $ \zeta_+^{\,-1} $  and the identifications  $ \,
\beta_{i,i} \cong t_{i,i} \, $,  $ \, \beta_{i,i+1} \cong t_{i,i+1} \, $.
In the other cases the result follows easily by induction on  $ \, j -
i \, $,  \, using the relations  $ \; \beta_{i,j} \; = \; \big( q -
q^{-1} \big)^{-1} \, \big( \, \beta_{i,k} \, \beta_{k,j} \, - \,
\beta_{k,j} \, \beta_{i,k} \big) \; \beta_{k,k}^{\;-1} \; $  (for
$ \, i < k < j \, $)  given in \S 2.6.
                                           \par
   Formulas (2.7) tell that the  $ \beta_{i,j} $'s  are quantum root
vectors too, for positive roots.  Similarly, for negative roots the
$ \gamma_{j,i} $'s  are involved.  Namely, the identifications
$ \zeta_-^{\,\pm 1} $  yield
  $$  \gamma_{i,i} \, = \, G_i^{\,-1} \; ,  \qquad
\gamma_{j,i} \, = \, - {(-q)}^{i-j} \, \fdot_{j,i}^- \, G_j^{\,-1}
\eqno \forall \;\; i < j \; \phantom{.} \quad \qquad (2.8)  $$
which are the analogues of (2.7).  Again this is proved by induction
on  $ \, j - i \, $:  \, the cases  $ \, j - i \leq 1 \, $  are
direct consequence of the description of  $ \zeta_-^{\,-1} $  and
the identifications  $ \, \gamma_{i,i} \cong t_{i,i} \, $,  $ \,
\gamma_{i+1,i} \cong t_{i+1,i} \, $,  \, while the inductive step
follow easily by means of the relations  $ \; \gamma_{j,i} \;
= \; \big( q - q^{-1} \big)^{-1} \, \big( \, \gamma_{k,i} \,
\gamma_{j,k} \, - \, \gamma_{j,k} \, \gamma_{k,i} \big) \;
\gamma_{k,k}^{\;+1} \; $  (for  $ \, j > k > i \, $)
given in \S 2.6.
                                           \par
   In order to compare (2.3) with (2.7) and (2.8) we must be able
to compare quantum root vectors with opposite superscripts.  The
tool is the unique  $ \kq $--algebra  antiautomorphism
  $$  \varPsi \, \colon \, \uqpgl \;{\buildrel \cong \over
{\lhook\joinrel\relbar\joinrel\relbar\joinrel\twoheadrightarrow}}\;
\uqpgl \; ,  \quad  E_i \mapsto E_i \; ,  \quad  F_i \mapsto F_i \; ,
\quad  G_j^{\,\pm 1} \mapsto G_j^{\,\mp 1}  \qquad \forall \;\, i, j  $$
which is clearly an involution; a straightforward computation shows that
  $$  \varPsi\big(E_{i,j}^\pm\big) \, = \, {(-q)}^{\mp(i-j+1)} \,
E_{i,j}^\mp \; ,  \qquad  \varPsi\big(F_{j,i}^\pm\big) \, = \,
{(-q)}^{\pm(i-j+1)} \, F_{j,i}^\mp  \qquad  \forall \;\, i < j \; .
\eqno (2.9)  $$
   \indent   Now comparing (2.3) with (2.7) and (2.8) using (2.9)
we get
  $$  \eqalignno{
   L_{i{}j}^+ \, = \, \varPsi\big(\gamma_{j,j}^{\,-1}
\, \gamma_{j,i} \, \gamma_{i,i}^{\,+1}\big)  \; ,  &  \qquad
L_{j\,i}^- \, = \, \varPsi\big(\beta_{i,i}^{\,+1} \,
\beta_{i,j} \, \beta_{j,j}^{\,-1}\big)  \hskip53pt  &
\qquad \quad \forall \;\, i \leq j \; ,   \qquad   (2.10)  \cr
   \gamma_{j,i} \, = \, \varPsi\big( (L_{ii}^+)^{-1}
\, L_{ij}^+ \, L_{jj}^+ \big)  \; ,  &  \qquad
\beta_{i,j} \, = \, \varPsi\big(L_{jj}^- \, L_{ji}^-
\, (L_{ii}^-)^{-1} \big)  \hskip53pt  &  \qquad
\quad \forall \;\, i \leq j \; .   \qquad   (2.11)  \cr }  $$

\vskip7pt

  {\bf 2.8 Presentation of  $ \, \utilde_q^{\scriptscriptstyle P}\!
(\gerg) \, $.}  \, Let again  $ \, G := {GL}_n \, $.  It is well known
that the  $ \kqqm $--integer  form  $ \fhat_q^{\scriptscriptstyle P}\!
\big[G\big] $  has the same presentation as  $ F_q^{\scriptscriptstyle P}
\!\big[G\big] $,  but over  $ \kqqm $  instead of  $ \kq \, $.  The
same holds similarly for  $ \fhat_q^{\scriptscriptstyle P}\!\big[
B_+ \big] $  and  $ \fhat_q^{\scriptscriptstyle P}\!\big[B_-\big] $.
In addition,  $ \, \fhat_q^{\scriptscriptstyle P}\!\big[B_\pm\big]
\cong \utilde_q^{\scriptscriptstyle P}\!(\gerb_\pm) \, $  and
$ \utilde_q^{\scriptscriptstyle P}\!(\gerg) $  is generated
by  $ \utilde_q^{\scriptscriptstyle P}\!(\gerb_+) $  and
$ \utilde_q^{\scriptscriptstyle P}\!(\gerb_-) $.
Therefore, the previous analysis implies that
$ \utilde_q^{\scriptscriptstyle P}\!(\gerg) $  as a
$ \kqqm $--algebra  is generated by the entries of the
$ q $--matrices  $ B $  and  $ \varGamma $  of Theorem
2.5.  The latter provides explicitly some relations (over
$ \kqqm $,  that is inside  $ \utilde_q^{\scriptscriptstyle P}
\!(\gerg) $  itself) among such generators, but these do  {\sl
not\/}  form a  {\sl complete\/}  set of relations: the general
mixed ones among  $ \beta_{i,j} $'s  and  $ \gamma_{r,s} $'s
are missing, as the ones in \S 2.6 do not make sense inside
$ \utilde_q^{\scriptscriptstyle P}\!(\gerg) \, $.  However,
since we know the relationship between these generators and
$ L $--operators  and we do know all relations among the
latter, we can eventually write down a complete set of
relations for the given generators!
   \hbox{This leads to the following presentation:}

\vskip7pt

\proclaim{Theorem 2.9}  {\sl (FRT-like presentation of
$ \, \utilde_q^{\scriptscriptstyle P}\!(\gergln) \, $)}
                                         \hfill\break
   \indent   $ \, \utilde_q^{\scriptscriptstyle P}\!(\gergln) \, $
is the unital  $ \kqqm $--algebra  with generators the entries of
the triangular matrices  $ \, B := \big( \beta_{i{}j} \big)_{i,j=1}^n
\, $  and  $ \, \varGamma := \big( \gamma_{i{}j} \big)_{i,j=1}^n \, $
and relations
  $$  \eqalignno{
   R \, B_2 \, B_1 \, = \, B_1 \, B_2 \, R \;\; ,  \hskip25pt  &
\hskip35pt  R \, \varGamma_2 \, \varGamma_1 \, = \, \varGamma_1
\, \varGamma_2 \, R  &  (2.12)  \cr
   R^{\text{op}} \, \varGamma^{\scriptscriptstyle D}_{\,1}
\, B^{\scriptscriptstyle D}_{\,2} \; = \;
B^{\scriptscriptstyle D}_{\,2} \,
\varGamma^{\scriptscriptstyle D}_{\,1} \, R^{\text{op}} \;\; ,
       \hskip9pt  &  \hskip21pt
  D_\beta \cdot D_\gamma \; = \; I \;= \; D_\gamma \cdot D_\beta
&   (2.13)  \cr }  $$
where  $ \, R := \sum_{i,j=1}^n \, q^{\delta_{i{}j}}
\, e_{i{}i} \otimes e_{j{}j} \, + \, \big( q - q^{-1} \big)
\, \textstyle{\sum_{1 \leq i < j \leq n}} \, e_{i{}j} \otimes
e_{j{}i} \, $,  $ \, X_1 := X \otimes I \, $,  $ \, X_2 := I
\otimes X \, $  (like in \S 2.1),  $ \, R^{\text{op}} :=
\sum_{i,j=1}^n \, q^{\delta_{i{}j}} \, e_{i{}i} \otimes e_{j{}j}
\, + \, \big( q - q^{-1} \big) \, \textstyle{\sum_{1 \leq i < j
\leq n}} \, e_{j{}i} \otimes e_{i{}j} \, $,  and  $ \, D_\beta
:= \text{\sl diag} \big( \beta_{1,1}, \dots, \beta_{n,n} \big) $,
$ \, D_\gamma := \text{\sl diag} \big( \gamma_{1,1}, \dots,
\gamma_{n,n} \big) $,  $ \, B^{\scriptscriptstyle D} := D_\beta^{+1}
\cdot B \cdot D_\beta^{-1} \, $,  $ \, \varGamma^{\scriptscriptstyle D}
:= D_\gamma^{-1} \cdot \varGamma \cdot D_\gamma^{+1} \, $.
                                          \par
   The first (compact) relation in (2.13) above is also equivalent to
  $$  {\textstyle \sum\limits_{i,k=1}^n} q^{\delta_{i,k}} \,
(e_{i,i} \otimes I) \Big( R^{\text{op}} \, \varGamma^-_{\,1} \,
B^+_{\,2} \Big) (I \otimes e_{k,k}) = \! {\textstyle \sum\limits_{j,
s = 1}^n} q^{\delta_{j,s}} \, (e_{j,j} \otimes I) \Big( B^-_{\,2}
\, \varGamma^+_{\,1} \, R^{\text{op}} \Big) (I \otimes e_{s,s})
\hskip11pt  \hfill (2.14)  $$
where  $ \, X^\pm := \big( q^{\pm \delta_{h,k}} \chi_{h,k} \big) \, $
for all  $ \, X \in \{ B, \varGamma \} $  (and  $ \, \chi \in \{ \beta,
\gamma \} $),  \, and in explicit, expanded form it is equivalent to
the set of relations (for all  $ \, i, k, j, s = 1, \dots, n $)
  $$  \hskip13pt   \hbox{ $ \eqalign{
   q^{\delta_{i,j}} \, \gamma_{i,k} \, \beta_{j,s} \;  &  + \;
\delta_{i>j} \, \big( q - q^{-1} \big) \, q^{\delta_{i,s} - \delta_{jk}}
\, \gamma_{j,k} \, \beta_{i,s} \;\; =  \cr
   &  \hskip25pt  = \;\; q^{\delta_{k,s}} \, \beta_{j,s} \, \gamma_{i,k}
\; + \; \delta_{s>k} \, \big( q - q^{-1} \big) \, q^{\delta_{i,s} -
\delta_{jk}} \beta_{j,k} \, \gamma_{i,s}  \cr } $ }   \eqno (2.15)  $$
where obviously  $ \, \delta_{h > k} := 1 \, $  if  $ \, h > k \, $
and  $ \, \delta_{h > k} := 0 \, $  if  $ \, h \not> k \, $.
                                   \par
   Furthermore,  $ \, \utilde_q^{\scriptscriptstyle P}\!
(\gergln) \, $  has the unique Hopf algebra structure given by
  $$  \Delta(X) \, = \, X \,\dot\otimes\, X \; ,  \qquad
\epsilon(X) \, = \, I \; ,  \qquad  S(X) = X^{-1}  \qquad
\;  \forall \;\, X \in \big\{ B, \varGamma \big\}  \; .
\eqno (2.16)  $$
\endproclaim

\demo{Proof} The commutation formulas in (2.12) and the Hopf formulas
in (2.16) are just the compact way to say that  $ B $  and  $ \varGamma $
are Hopf  $ q $--matrices.  The second half of (2.13) instead is nothing
but another way of writing (2.2).
                                      \par
   Moreover, the first half of (2.13) arises from the similar compact
relation for  $ L $--operators  and the link between the latter and
the present generators.  Indeed, merging (2.10) in the last identity
in (2.4) we get
  $$  R \cdot \varPsi \big( D_\gamma^{\,-1}
\varGamma^{\scriptscriptstyle T} D_\gamma^{\,+1} \big)_1
\cdot \varPsi \big( D_\beta^{\,+1} B^{\scriptscriptstyle T}
D_\beta^{\,-1} \big)_2 \; = \;
\varPsi \big( D_\beta^{\,+1} B^{\scriptscriptstyle T}
D_\beta^{\,-1} \big)_2 \cdot \varPsi \big( D_\gamma^{\,-1}
\varGamma^{\scriptscriptstyle T} D_\gamma^{\,+1} \big)_1
\cdot R  $$
(where a superscript  $ T $  means ``transpose'').  Using the
fact that  $ \varPsi $  is an algebra antiautomorphism and
extending its action to  $ \, \varPsi(R) = R \, $  we then
argue
  $$  \varPsi \Big( \hskip-2pt \big( D_\beta^{\,+1} B \, D_\beta^{\,-1}
\big)_2 \cdot \big( D_\gamma^{\,-1} \varGamma \, D_\gamma^{\,+1} \big)_1
\cdot R^{\text{op}} \Big) \; = \;
   \varPsi \Big( R^{\text{op}} \cdot \big( D_\gamma^{\,-1} \varGamma
\, D_\gamma^{\,+1} \big)_1 \cdot \big( D_\beta^{\,+1} B \,
D_\beta^{\,-1} \big)_2 \Big)  $$
whence eventually (2.13) follows at once because  $ \; \varPsi^2 =
\text{\sl id} \; $.
                                      \par
   Finally, expanding (2.13) one gets explicitly (for all  $ \, i, k, j,
s = 1, \dots, n $)
  $$  \eqalign{
   q^{\delta_{i,j}} \,  &  \gamma_{i,i}^{\,-1} \, \gamma_{i,k}
\, \gamma_{k,k}^{\,+1} \, \beta_{j,j}^{\,+1} \, \beta_{j,s} \,
\beta_{s,s}^{\,-1} \, + \, \delta_{i>j} \, \big( q - q^{-1} \big)
\, \gamma_{j,j}^{\,-1} \, \gamma_{j,k} \, \gamma_{k,k}^{\,+1} \,
\beta_{i,i}^{\,+1} \, \beta_{i,s} \, \beta_{s,s}^{\,-1} \; =  \cr
   &  \hskip15pt  = \; q^{\delta_{k,s}} \, \beta_{j,j}^{\,+1} \,
\beta_{j,s} \, \beta_{s,s}^{\,-1} \, \gamma_{i,i}^{\,-1} \, \gamma_{i,k}
\, \gamma_{k,k}^{\,+1} \, + \, \delta_{s>k} \, \big( q - q^{-1} \big)
\, \beta_{j,j}^{\,+1} \, \beta_{j,k} \, \beta_{k,k}^{\,-1} \,
\gamma_{i,i}^{\,-1} \, \gamma_{i,s} \, \gamma_{s,s}^{\,+1} \; .  \cr }  $$
From this, making repeated use of all the relations encoded in (2.12) and
in the second half of (2.13) one can cancel out all ``diagonal'' factors,
i.e.~those of type  $ \beta_{\ell,\ell} $  or  $ \gamma_{\ell,\ell} \, $.
The outcome is (for all  $ \, i, k, j, s = 1, \dots, n $)
  $$  \eqalign{
   q^{\delta_{i,j}} \, \gamma_{i,k} \, \beta_{j,s} \;  &  + \;
\delta_{i>j} \, \big( q - q^{-1} \big) \, q^{\delta_{i,s} - \delta_{jk}}
\, \gamma_{j,k} \, \beta_{i,s} \;\; =  \cr
   &  \hskip75pt  = \;\; q^{\delta_{k,s}} \, \beta_{j,s} \, \gamma_{i,k}
\; + \; \delta_{s>k} \, \big( q - q^{-1} \big) \, q^{\delta_{i,s} -
\delta_{jk}} \beta_{j,k} \, \gamma_{i,s}  \cr }  $$
that is exactly the set of relations (2.15).  As a last step,
manipulating a bit the exponents of  $ q $  one gets  (for
$ \, i, k, j, s = 1, \dots, n $)
  $$  \hskip-3pt   \hbox{ $ \eqalign{
   &  q^{2 \, \delta_{i,k}} \Big( q^{\delta_{i,j}}
\, \big( q^{-\delta_{i,k}} \, \gamma_{i,k}\big) \,
\big( q^{+\delta_{j,s}} \beta_{j,s}\big) \, + \,
\delta_{i>j} \, \big( q - q^{-1} \big) \,
\big( q^{-\delta_{j,k}} \, \gamma_{j,k} \big)
\, \big( q^{+\delta_{i,s}} \beta_{i,s} \big) \Big) \, =  \cr
   &  \;\; = \; q^{2 \, \delta_{j,s}} \Big( q^{\delta_{k,s}} \,
\big( q^{-\delta_{j,s}} \beta_{j,s} \big) \, \big( q^{+\delta_{i,k}}
\, \gamma_{i,k} \big) \, + \, \delta_{s>k} \, \big( q - q^{-1} \big) \,
\big( q^{-\delta_{j,k}} \beta_{j,k} \big) \, \big( q^{+\delta_{i,s}} \,
\gamma_{i,s}\big) \Big)  \cr } $ }   \hskip-2pt  (2.15)  $$
when written in compact form yields exactly (2.14).   \qed
\enddemo

\vskip7pt

   {\it  $ \underline{\hbox{\it Remark}} $:} \; the argument used to
argue formulas (2.13) from the last identity in (2.4) may be also
applied to the first two identities therein.  This yields relations
among the  $ \beta_{ij} $'s  and among the  $ \gamma_{ji} $'s  which
are different from, but equivalent to, formulas (2.12).

\vskip7pt

\proclaim{Corollary 2.10} \, The Poisson-Hopf\/  $ \Bbbk $--algebra
$ \, \utilde_1^{\scriptscriptstyle P}\!(\gergln) $  is the polynomial,
Laurent-polyno\-mial algebra in the variables  $ \, \big\{
\overline{\beta}_{i,j} \big\}_{1 \leq i \leq j \leq n} \bigcup
\big\{ \overline{\gamma}_{j,i} \big\}_{1 \leq i \leq j \leq n} \; $,
\; the  $ \beta_{\ell{}\ell} $'s  and  $ \gamma_{ii} $'s  being
invertible, with relations  $ \; \beta_{ii}^{\,\pm 1} \! =
\gamma_{ii}^{\,\mp 1} \, (\, \forall \, i ) $,  Hopf structure
     \hbox{given (in compact notation) by}
  $$  \Delta\big(\,\overline{X}\,\big) \, = \, \overline{X} \,\dot\otimes\,
\overline{X} \;\, ,  \qquad  \epsilon\big(\,\overline{X}\,\big) \, = \, I
\;\, ,  \qquad  S\big(\,\overline{X}\,\big) = \overline{X}^{\;-1}  \qquad
\;  \forall \;\, X \in \big\{ B, \varGamma \big\}  $$
(with  $ B $  and  $ \varGamma $  as in Theorem 2.9) and with the unique
Poisson structure such that
  $$  \hbox{ $ \eqalign{
%
%
  \big\{ \overline{x}_{i,h} \, , \overline{x}_{i,\ell} \big\}
= \, \overline{x}_{i,h} \, \overline{x}_{i,\ell} \; ,  \hskip8pt
\big\{ \overline{x}_{h,j} \, , \overline{x}_{\ell,j} \big\} =  &
\,\, \overline{x}_{h,j} \, \overline{x}_{\ell,j} \; ,  \hskip8pt
\big\{ \overline{x}_{h,h} \, , \overline{x}_{\ell,\ell} \big\}
= \, 0   \hskip12pt  (\, h < \ell \,)
 \cr
   \big\{ \overline{x}_{i,j} , \overline{x}_{h,k} \big\}
= \, 0   \hskip11pt (\, i < h \, , j > k \,) \; ,  \qquad  &
\big\{ \overline{x}_{i,j} \, , \overline{x}_{h,k} \big\} =
2 \; \overline{x}_{i,k} \, \overline{x}_{h,j}  \hskip11pt
(\, i < h \, , j < k \,)  \cr } $ }   \hskip2pt (2.17)   $$
with either all  $ x_{p\,q} $'s  being  $ \beta_{p\,q} $'s  (and
$ \, \beta_{p\,q} := 0 \, $  for all  $ \, p > q \, $)  or all
$ x_{p\,q} $'s  being  $ \gamma_{p\,q} $'s  (and  $ \, \gamma_{p\,q}
:= 0 \, $  for all  $ \, p < q \, $),  and
  $$  \big\{ \overline{\beta}_{j,s} \, , \overline{\gamma}_{i,k} \big\}
\; = \; (\delta_{i,j} - \delta_{k,s}) \cdot \overline{\beta}_{j,s}
\, \overline{\gamma}_{i,k} \, + \, 2 \; \delta_{i>j} \cdot
\overline{\gamma}_{j,k} \, \overline{\beta}_{i,s} \, -
\, 2 \; \delta_{s>k} \cdot \overline{\beta}_{j,k}
\, \overline{\gamma}_{i,s} \quad .   \eqno (2.18)  $$
%
%
   In particular  $ \, \utilde_1^{\scriptscriptstyle P}\!
(\gergln) \cong F\big[({GL}_n)^*_{\scriptscriptstyle P}\big] \, $  as
Poisson Hopf algebras, where  $ ({GL}_n)^*_{\scriptscriptstyle P} $
is the algebraic group of pairs of matrices  $ \big( \varGamma, B \big) $
where  $ \varGamma $,  resp.~$ B $,  is a lower triangular, resp.~upper
triangular, invertible matrix, and the diagonals of  $ \, \varGamma $
and  $ B $  are inverse to each other, with the Poisson structure dual
to the Lie bialgebra of  $ \, \gergln \, $.
\endproclaim

\demo{Proof} If we write  $ \, \overline{x} := x \mod (q-1)
\, \utilde_q^{\scriptscriptstyle P}\!(\gergln) \, $  for
every  $ \, x \in \utilde_q^{\scriptscriptstyle P}\!(\gergln)
\, $,  \, then setting  $ \, q = 1 \, $  in the presentation
of $ \utilde_q^{\scriptscriptstyle P}\!(\gergln) $
of Theorem 2.9 yields a presentation for  $ \,
\utilde_1^{\scriptscriptstyle P}\!(\gergln) \, $.
The latter is a commutative, polynomial
Laurent-polynomial algebra as claimed, whence
$ \, \utilde_1^{\scriptscriptstyle P}\!(\gergln) \cong
F\big[({GL}_n)^*_{\scriptscriptstyle P}\big] \, $  as
algebras, via an isomorphism which for all  $ \, i \leq j \, $
maps  $ \, \overline{\beta_{ij}} := \beta_{ij} \mod (q-1) \,
\utilde_q^{\scriptscriptstyle P}\!(\gergln) \, $  to the matrix
coefficient corresponding to the  $ (i,j) $--th  entry of the
matrix  $ B $  in a pair  $ (\varGamma,B) $  as in the claim,
and maps  $ \, \overline{\gamma_{ji}} := \gamma_{ji} \mod (q-1) \,
\utilde_q^{\scriptscriptstyle P}\!(\gergln) \, $  to the matrix
coefficient corresponding to the  $ (j,i) $--th  entry of the matrix
$ \varGamma $  in a pair  $ (\varGamma,B) $.  The formulas for the
Hopf structure in  $ \, \utilde_q^{\scriptscriptstyle P}\!(\gergln)
\, $  imply that this is also an isomorphism of Hopf algebras, for
the Hopf structure on right hand side induced by the group structure
of  $ ({GL}_n)^*_{\scriptscriptstyle P} \, $.
                                                    \par
   Since  $ \utilde_1^{\scriptscriptstyle P}\!(\gergln) $  is
commutative, it inherits from  $ \utilde_q^{\scriptscriptstyle P}\!
(\gergln) $  the unique Poisson bracket given by the rule  $ \;
\displaystyle{ \big\{ \overline{x} \, , \overline{y} \,\big\}
:= {{\, x \, y - y \, x \,} \over {\, q - 1 \,}} \mod (q-1) \,
\utilde_q^{\scriptscriptstyle P}\!(\gergln) } \phantom{\Big|} $,
\; for all  $ \, x, y \in \utilde_q^{\scriptscriptstyle P}\!
(\gergln) \, $.  Then the Poisson brackets in (2.18) come directly
from (2.15), while all those in (2.17) spring out of the commutation
formulas among the  $ \beta_{ij} $'s  and among the  $ \gamma_{ji} $'s
in (2.11).
                                                    \par
   Finally, checking that this Poisson structure on the algebraic group
$ ({GL}_n)^*_{\scriptscriptstyle P} $  is exactly the one dual to the
Lie bialgebra structure of  $ \gergl_{\,n} $  is just a matter of
bookkeeping.   \qed
\enddemo

\vskip7pt

  {\bf 2.11 The quantum Frobenius morphisms  $ \,
F \big[ ({GL}_n)^*_{\scriptscriptstyle P} \big] \cong
\, \utilde_1^{\scriptscriptstyle P}\!(\gergln) \,
\lhook\joinrel\relbar\joinrel\relbar\joinrel\relbar\joinrel\rightarrow
\, \utilde_\varepsilon^{\scriptscriptstyle P}\!(\gergln)
\, $.}  \, Let  $ \, \keps \, $  be the extension of  $ \Bbbk $  by a
primitive  $ \ell $--th  root of 1, say  $ \varepsilon \, $.  Since
$ \, \utilde_q^{\scriptscriptstyle P}\!(\gergln) \, $
is generated by copies of  $ \, \utilde_q^{\scriptscriptstyle P}\!
(\gerb_+) \cong \fhat_q^{\scriptscriptstyle P}\!\big[B_+\big] \, $
and  $ \, \utilde_q^{\scriptscriptstyle P}\!(\gerb_-) \cong
\fhat_q^{\scriptscriptstyle P}\!\big[B_-\big] \, $,
\, taking specializations the same is true for
$ \utilde_\varepsilon^{\scriptscriptstyle P}\!
(\gergln) \, $;  \, in particular the latter
is presented like in Theorem 2.9 but with
$ \, q = \varepsilon \, $.
                                            \par
   In addition, the quantum Frobenius morphisms  $ \, F\big[ {GL}_n
\big] \cong \fhat_1^{\scriptscriptstyle P}\!\big[{GL}_n\big]
\lhook\joinrel\relbar\joinrel\relbar\joinrel\rightarrow
\fhat_\varepsilon^{\scriptscriptstyle P}\!\big[{GL}_n\big] \, $  and
$ \, F\big[B_\pm\big] \cong \fhat_1^{\scriptscriptstyle P}\! \big[
B_\pm \big] \lhook\joinrel\relbar\joinrel\relbar\joinrel\rightarrow
\fhat_\varepsilon^{\scriptscriptstyle P}\!\big[B_\pm\big] \, $
have a pretty neat description, as they are given by  $ \, t_{i,j}
\mapsto t_{i,j}^{\;\ell} \, $  (hereafter we denote by the same
symbol an element in a quantum algebra and its corresponding coset
after any specialization); see, for instance, [PW] for details.
As we mentioned in \S 1.6, the morphism  $ \, F \big[
({GL}_n)^*_{\scriptscriptstyle P} \big] \cong
\, \utilde_1^{\scriptscriptstyle P}\!(\gergln) \,
\lhook\joinrel\relbar\joinrel\relbar\joinrel\relbar\joinrel\rightarrow
\, \utilde_\varepsilon^{\scriptscriptstyle P}\!(\gergln)
\, $  is determined by its restriction to the quantum Borel subalgebras,
hence to the copies of  $ \fhat_1^{\scriptscriptstyle P}\! \big[ B_+
\big] $  and  $ \fhat_1^{\scriptscriptstyle P}\! \big[ B_- \big] $
which generate  $ \utilde_1^{\scriptscriptstyle P}\!(\gergln) $.
When reformulated in light of Theorem 2.9, this implies

\vskip7pt

\proclaim{Theorem 2.12} The quantum Frobenius morphism
   \hbox{$ \; F \big[ ({GL}_n)^*_{\scriptscriptstyle P} \big] \cong
\, \utilde_1^{\scriptscriptstyle P}\!(\gergln)
\lhook\joinrel\relbar\joinrel\relbar\joinrel\rightarrow
\utilde_\varepsilon^{\scriptscriptstyle P}\!(\gergln) $}
\; is given by  $ \; \beta_{i,j} \mapsto \beta_{i,j}^{\;\ell} \; $,
$ \; \gamma_{j,i} \mapsto \gamma_{j,i}^{\;\ell} \; $,  \; for all
$ \, i \leq j \; $.   \qed
\endproclaim

\vskip1,3truecm

\centerline {\bf \S \; 3 \ The case of  $ \gersln \, $ }

\vskip10pt

  {\bf 3.1 From  $ \gergln $  to  $ \gersln \, $.} \, In this section,
we consider  $ \, \gerg = \gersln \, $  and  $ \, G = {SL}_n \, $.
The constructions and results of \S 2 about  $ \gergln $  essentially
duplicate into the like for  $ \gersln \, $,  \, up to minor details.
In this section we shall draw these results,
      \hbox{shortly explaining the changes in order.}
                                      \par
   First, the ideal generated by  $ \, (L_n - 1) \, $  in  $ \uqpgl $
is a  {\sl Hopf ideal\/}:  \, then we define  $ \uqpsl $  as the
quotient Hopf  $ \kq $--algebra  $ \, \uqpsl := \uqpgl \Big/ \big(
L_n - 1 \big) \, $.  With like notation (see \S 2.3) for generators
of  $ \uqpgl $  and their images in  $ \uqpsl $,  we define  $ \uqqsl $
as the  $ \kq $--subalgebra  of  $ \uqpsl $  generated by  $ \, \big\{
F_i, K_i^{\pm 1}, E_i \big\}_{i = 1, \dots, n-1} \, $;  \, this is
also the image of  $ \uqqgl $  when mapping  $ \uqpgl $  onto
$ \uqpsl $.  In this setting,  $ \uqpbp $,  resp.~$ \uqpbm $,
is the  $ \kq $--subalgebra  of  $ \uqpsl $  generated by
$ \big\{ L_i^{\pm 1}, E_i \big\}_{i = 1, \dots, n-1} \, $,
resp.~$ \, \big\{ F_i, L_i^{\pm 1} \big\}_{i = 1, \dots, n-1}
\, $,  \, whereas  $ \uqqbp $,  resp.~$ \uqqbm $,  instead is
the  $ \kq $--subalgebra of
      $ \uqqsl $  generated by\break
$ \, \big\{ K_i^{\pm 1} , E_i \big\}_{i = 1, \dots, n-1} \, $,
resp.~$ \, \big\{ F_i, K_i \big\}_{i = 1, \dots, n-1} \, $.  All
these are  {\sl Hopf\/}  subalgebras of  $ \uqpsl $  and  $ \uqqsl $,
and Hopf algebra quotients of the similar quantum Borel algebras for
$ \gergl_n \, $.
                                          \par
   In this context, we can repeat step by step the construction made
for  $ \, \gergl_n \, $,  up to minimal details (namely, taking into
account everywhere the relation  $ \, L_n = 1 \, $);  in particular,
in quantum function algebras the additional relation  $ \, t_{1,1}
\, t_{2,2} \cdots t_{n,n} = 1 \, $  has to be taken into account.
Otherwise, the results for the  $ \gersl_n $  case can be immediately
argued from the corresponding results for  $ \gergl_n \, $.  The first
of these results   --- analogue to Theorem 2.5 ---  is

\vskip7pt

\proclaim{Theorem 3.2}  {\sl (``short'' FRT-like presentation of
$ \uqpsl \, $)}
                                         \hfill\break
   \indent   $ \, \uqpsl $  is the quotient algebra of  $ \,
\uqpgl $  modulo the two-sided ideal  $ I $  generated by  $ \, \Big(
\prod\nolimits_{i=1}^n \beta_{ii} - 1 \Big) \, $  (or by  $ \, \Big(
\prod\nolimits_{j=1}^n \gamma_{jj} - 1 \big) \, $,  \, which gives the
same).  Moreover,  $ I $  is a Hopf ideal of  $ \, \uqpgl $,  therefore
$ \, \uqpsl $  inherits from  $ \uqpgl $  a structure of quotient Hopf
algebra, given by formulas like in Theorem 2.5 (with the obvious,
additional simplifications).  In particular,  $ \, \uqpsl $  has the
same presentation as  $ \, \uqpgl $  in Theorem 2.5  {\sl plus}  the
additional relation  $ \; \beta_{1,1} \, \beta_{2,2} \cdots \beta_{n,n}
= 1 \, $,  \, or  $ \; \gamma_{1,1} \, \gamma_{2,2} \cdots \gamma_{n,n}
= 1 \, $.   \qed
\endproclaim

%
%
\vskip3pt

  {\bf 3.3 Quantum root vectors,  $ L $--operators  and new
generators for  $ \utilde_q^{\scriptscriptstyle P}\!(\gersln) $.}
\, Definitions imply that the Hopf algebra epimorphism  $ \, \uqpgl
\relbar\joinrel\twoheadrightarrow \uqpsl \, $  maps any quantum root
vector   --- say  $ E_{i,j} $  or  $ F_{j,i} $  ---   in  $ \uqpgl $
onto a corresponding quantum root vector in  $ \uqpsl \, $,  for
which we use the like notation.  A similar result clearly holds
for each  $ L $--operator   --- in  $ \uqpgl $  ---   too, whose
image in  $ \uqpsl $  we still denote by the same symbol.  The
discussion in \S\S 2.7--9 can then be repeated  {\it verbatim},
in particular formulas (2.3) through (2.11) hold true within
$ \uqpsl $  as well.  The outcome then is the analogue of
Theorem 2.9   --- which can also be deduced directly from the
latter, since  $ \utilde_q^{\scriptscriptstyle P}\!(\gergln) $
maps onto  $ \utilde_q^{\scriptscriptstyle P}\!(\gersln) $  ---
and its immediate corollary, namely

\vskip7pt

\proclaim{Theorem 3.4}  {\sl (FRT-like presentation of
$ \, \utilde_q^{\scriptscriptstyle P}\!(\gersln) \, $)}
                                         \hfill\break
   \indent   $ \, \utilde_q^{\scriptscriptstyle P}\!(\gersln) \, $
is the unital  $ \kqqm $--algebra  with generators the entries of
the triangular matrices  $ \, B := \big( \beta_{i{}j} \big)_{i,j=1}^n
\, $  and  $ \, \varGamma := \big( \gamma_{i{}j} \big)_{i,j=1}^n \, $
and relations (notations as in Theorem 2.9)
  $$  \eqalignno{
   R \, B_2 \, B_1 \, = \, B_1 \, B_2 \, R \;\; ,  \hskip25pt  &
\hskip35pt  R \, \varGamma_2 \, \varGamma_1 \, = \, \varGamma_1
\, \varGamma_2 \, R  &  (3.1)  \cr
   R^{\text{op}} \, \varGamma^{\scriptscriptstyle D}_{\,1}
\, B^{\scriptscriptstyle D}_{\,2} \; = \;
B^{\scriptscriptstyle D}_{\,2} \,
\varGamma^{\scriptscriptstyle D}_{\,1} \, R^{\text{op}} \;\; ,
       \hskip9pt  &  \hskip21pt
  D_\beta \cdot D_\gamma \; = \; I \;= \; D_\gamma \cdot D_\beta
&   (3.2)  \cr
   \text{\sl det}\,(D_\beta) \; = \;  &
       \; 1 \; = \; \text{\sl det}\,(D_\gamma)   &  (3.3)  \cr }  $$
   \indent   The first (compact) relation in (2.13) above is
equivalent to
  $$  {\textstyle \sum\limits_{i,k=1}^n} q^{\delta_{i,k}} \,
(e_{i,i} \otimes I) \Big( R^{\text{op}} \, \varGamma^-_{\,1} \,
B^+_{\,2} \Big) (I \otimes e_{k,k}) = \! {\textstyle \sum\limits_{j,
s = 1}^n} q^{\delta_{j,s}} \, (e_{j,j} \otimes I) \Big( B^-_{\,2}
\, \varGamma^+_{\,1} \, R^{\text{op}} \Big) (I \otimes e_{s,s})
\hskip11pt  \hfill (3.4)  $$
and in expanded form it is equivalent to the set of relations
(for all  $ \, i, k, j, s = 1, \dots, n $)
  $$  \hskip13pt   \hbox{ $ \eqalign{
   q^{\delta_{i,j}} \, \gamma_{i,k} \, \beta_{j,s} \;  &  + \;
\delta_{i>j} \, \big( q - q^{-1} \big) \, q^{\delta_{i,s} - \delta_{jk}}
\, \gamma_{j,k} \, \beta_{i,s} \;\; =  \cr
   &  \hskip25pt  = \;\; q^{\delta_{k,s}} \, \beta_{j,s} \, \gamma_{i,k}
\; + \; \delta_{s>k} \, \big( q - q^{-1} \big) \, q^{\delta_{i,s} -
\delta_{jk}} \beta_{j,k} \, \gamma_{i,s}  \quad .  \cr } $ }
\eqno (3.5)  $$
   \indent   Furthermore,  $ \, \utilde_q^{\scriptscriptstyle P}\!
(\gersln) \, $  has the unique Hopf algebra structure given by
  $$  \Delta(X) \, = \, X \,\dot\otimes\, X \; ,  \qquad
\epsilon(X) \, = \, I \; ,  \qquad  S(X) = X^{-1}  \qquad
\;  \forall \;\, X \in \big\{ B, \varGamma \big\}  \; .
\qed   \eqno (3.6)  $$
\endproclaim

\vskip7pt

\proclaim{Corollary 3.5} \, The Poisson-Hopf\/  $ \Bbbk $--algebra
$ \, \utilde_1^{\scriptscriptstyle P}\!(\gersln) $  is the polynomial
algebra in the variables  $ \, \big\{ \overline{\beta}_{i,j} \big\}_{1
\leq i \leq j \leq n} \bigcup \big\{ \overline{\gamma}_{j,i} \big\}_{1
\leq i \leq j \leq n} \; $  modulo the relations  $ \;
\overline{\beta}_{1,1} \, \overline{\beta}_{2,2} \cdots
\overline{\beta}_{n,n} \, = \, 1 \, $,  $ \; \overline{\gamma}_{1,1}
\, \overline{\gamma}_{2,2} \cdots \overline{\gamma}_{n,n} = \, 1 \, $,
$ \; \overline{\beta}_{i,i} \, \overline{\gamma}_{i,1i} = \, 1 \; $
(for all  $ \, i = 1, \dots n \, $),  \, Hopf structure given by
  $$  \Delta\big(\,\overline{X}\,\big) \, = \, \overline{X} \,\dot\otimes\,
\overline{X} \;\, ,  \qquad  \epsilon\big(\,\overline{X}\,\big) \, = \, I
\;\, ,  \qquad  S\big(\,\overline{X}\,\big) = \overline{X}^{\;-1}  \qquad
\;  \forall \;\, X \in \big\{ B, \varGamma \big\}  $$
(with  $ B $  and  $ \varGamma $  as in Theorem 3.4) and with the unique
Poisson structure such that
  $$  \hbox{ $ \eqalign{
  \big\{ \overline{x}_{i,h} \, , \overline{x}_{i,\ell} \big\}
= \, \overline{x}_{i,h} \, \overline{x}_{i,\ell} \; ,  \hskip8pt
\big\{ \overline{x}_{h,j} \, , \overline{x}_{\ell,j} \big\} =  &
\,\, \overline{x}_{h,j} \, \overline{x}_{\ell,j} \; ,  \hskip8pt
\big\{ \overline{x}_{h,h} \, , \overline{x}_{\ell,\ell} \big\}
= \, 0   \hskip12pt  (\, h < \ell \,)  \cr
   \big\{ \overline{x}_{i,j} , \overline{x}_{h,k} \big\}
= \, 0   \hskip11pt (\, i < h \, , j > k \,) \; ,  \qquad  &
\big\{ \overline{x}_{i,j} \, , \overline{x}_{h,k} \big\} =
2 \; \overline{x}_{i,k} \, \overline{x}_{h,j}  \hskip11pt
(\, i < h \, , j < k \,)  \cr } $ }   \hskip2pt (3.7)   $$
with either all  $ x_{p\,q} $'s  being  $ \beta_{p\,q} $'s  (and
$ \, \beta_{p\,q} := 0 \, $  for all  $ \, p > q \, $)  or all
$ x_{p\,q} $'s  being  $ \gamma_{p\,q} $'s  (and  $ \, \gamma_{p\,q}
:= 0 \, $  for all  $ \, p < q \, $),  and
  $$  \big\{ \overline{\beta}_{j,s} \, , \overline{\gamma}_{i,k} \big\}
\; = \; (\delta_{i,j} - \delta_{k,s}) \cdot \overline{\beta}_{j,s}
\, \overline{\gamma}_{i,k} \, + \, 2 \; \delta_{i>j} \cdot
\overline{\gamma}_{j,k} \, \overline{\beta}_{i,s} \, -
\, 2 \; \delta_{s>k} \cdot \overline{\beta}_{j,k}
\, \overline{\gamma}_{i,s} \quad .   \eqno (3.8)  $$
   In particular  $ \, \utilde_1^{\scriptscriptstyle P}\!
(\gersln) \cong F\big[({SL}_n)^*_{\scriptscriptstyle P}\big] \, $  as
Poisson Hopf algebras, where  $ ({SL}_n)^*_{\scriptscriptstyle P} $
is the algebraic group of pairs of matrices  $ \big( \varGamma, B \big) $
where  $ \varGamma $,  resp.~$ B $,  is a lower, resp.~upper, triangular
matrix with determinant equal to 1, and the diagonals of  $ \, \varGamma $
and  $ B $  are inverse to each other, with the Poisson structure dual to
the Lie bialgebra of  $ \, \gersln \, $.   \qed
\endproclaim

\vskip7pt

  {\bf 3.7 The quantum Frobenius morphisms  $ \,
F \big[ ({SL}_n)^*_{\scriptscriptstyle P} \big] \cong
\, \utilde_1^{\scriptscriptstyle P}\!(\gersln) \,
\lhook\joinrel\relbar\joinrel\relbar\joinrel\relbar\joinrel\rightarrow
\, \utilde_\varepsilon^{\scriptscriptstyle P}\!(\gersln) \, $.}
\, Once again, for quantum Frobenius morphisms one can repeat  {\it
verbatim\/}  the discussion made for  $ \uqpgl $  for the case of
$ \uqpsl $  as well, via minimal changes where needed.  Otherwise,
the results in the  $ \gergl_n $  case induce similar results in
the  $ \gersl_n $  case via the defining epimorphism  $ \, \uqpgl
\relbar\joinrel\twoheadrightarrow \uqpsl \, $.  Indeed, the latter
is clearly compatible (in the obvious sense) with specializations
at roots of 1; therefore, the specializations of the epimorphism
itself yield the following commutative diagram
  $$  \CD
   F \big[ ({GL}_n)^*_{\scriptscriptstyle P} \big] \cong
\, \utilde_1^{\scriptscriptstyle P}\!(\gergln)  @>>>
\utilde_\varepsilon^{\scriptscriptstyle P}\!(\gergln)  \\
   @VVV   @VVV  \\
   F \big[ ({SL}_n)^*_{\scriptscriptstyle P} \big] \cong
\, \utilde_1^{\scriptscriptstyle P}\!(\gersln)  @>>>
\utilde_\varepsilon^{\scriptscriptstyle P}\!(\gersln)
      \endCD  $$
(for  $ \varepsilon $  any root of 1) in which the vertical arrows
are the above mentioned specialized epimorphisms and the horizontal
ones are the quantum Frobenius (mono)morphisms.
                                            \par
   This yields at once the following analogue of Theorem 2.12:

\vskip7pt

\proclaim{Theorem 3.8} The quantum Frobenius morphism
   \hbox{$ \; F \big[ ({SL}_n)^*_{\scriptscriptstyle P} \big] \cong
\, \utilde_1^{\scriptscriptstyle P}\!(\gersln)
\lhook\joinrel\relbar\joinrel\relbar\joinrel\rightarrow
\utilde_\varepsilon^{\scriptscriptstyle P}\!(\gersln) $}
\; is given by  $ \; \beta_{i,j} \mapsto \beta_{i,j}^{\;\ell} \; $,
$ \; \gamma_{j,i} \mapsto \gamma_{j,i}^{\;\ell} \; $,  \; for all
$ \, i \leq j \; $.   \qed
\endproclaim

\vskip1,5truecm

\vskip39pt

\enddocument